\documentclass[11pt]{amsart}
\input epsf
\usepackage{subfigure,graphicx, caption,psfrag} %pst-node}

\usepackage[percent]{overpic} % to put overlay latex on images -- 
\title{Circle Packing with Generalized Branching} 

\author{James Ashe, Edward Crane and Kenneth Stephenson}
\date{\today}

\newtheorem{theorem}{Theorem}
\newtheorem{lemma}[theorem]{Lemma}

\theoremstyle{definition}

\theoremstyle{remark}

\newcommand{\fB}{\mathbf B} % classical Blaschke product
\newcommand{\fb}{\mathbf b} % discrete Blaschke product
\newcommand{\fw}{\mathbf w} % classical Weierstrass 
\newcommand{\fW}{\mathbf W} % discrete Weierstrass 
\newcommand{\fA}{\mathbf A} % classical Ahlfors
\newcommand{\fa}{\mathbf a} % discrete Ahlfors
\newcommand{\ff}{\mathbf f} % generic disc analytic
\newcommand{\fF}{\mathbf F} % generic classical analytic

\newcommand{\F}[1]{Figure~\ref{F:#1}}    % Figure reference

\newcommand{\ntn}[1]{}   % fake for when I don't want these to print.
\newcommand{\hide}[1]{}	    % for hiding: add #1 if you want these to show

\newcommand{\bbb}{\mathbb}      % black board bold
\newcommand{\CP}{\text{\tt CirclePack}} % software name

\newcommand{\bC}{\bbb C}         % complex plane
\newcommand{\bD}{\bbb D}         % unit disc (hyperbolic plane)
\newcommand{\bP}{\bbb P}         % Riemann sphere
\newcommand{\wK}{\widetilde K}

\newcommand{\wP}{\widetilde{P}}

\newcommand{\wR}{\widetilde{R}}

\begin{document}
\maketitle

\begin{center}{\it Dedicated to C. David Minda on the occasion of 
his retirement}
\end{center}

Classical analytic function theory is at the heart of David Minda's
research and of many of the results in this volume. It has been a
pleasure in recent years to find that simple patterns of circles
called circle packings could find themselves in such tight company
with this classical theory. David himself contributed to this topic in
\cite{MR91} and on his retirement will surely have time to dive back
into it.

Let us briefly review the circle packing story line. It began with Bill
Thurston's observation that for every abstract triangulation $K$ of a
topological sphere there exists an essentially unique configuration of
circles with mutually disjoint interiors on the Riemann sphere $\bP$
whose pattern of tangencies is encoded in $K$. That is, there's a
circle packing $P$ for $K$ in $\bP$. Based on this rigidity and his
intuition, Thurston made a remarkable proposal at the 1985 Conference
in Celebration of de~Branges' Proof of the Bieberbach Conjecture:
namely, that one could use such circle packings to approximate
conformal mappings.  The subsequent proof of his conjecture by Burt
Rodin and Dennis Sullivan \cite{RS87} established circle packing as a
topic and opened its most widely known aspect, the approximation of
classical analytic functions.

As this approximation theory developed, a second aspect that we will
call {\sl discrete analytic function theory}, began to emerge. For it
became increasingly clear that classical phenomena were already at
play within circle packing --- mappings between circle packings not
only approximated analytic functions, they also mimicked them. The
literature shows an ever growing list of conformal notions being
realized discretely and often with remarkable geometric fidelity:
moving circle packing into the hyperbolic geometry of $\bD$ led to
infinite packings and the consequent classical type conditions --- the
spherical, hyperbolic, and euclidean trichotomy --- and from that
came the discrete uniformization theorem, discrete Riemann surfaces
and covering theory, random walks, and then notions of branch points
and boundary conditions allowed for discrete versions of familiar
classes of functions, polynomials, exponentials, and the Blaschke
products, Ahlfors functions, and Weierstrass functions that play
their roles in this paper.

Part and parcel in these developments has been a third aspect,
computation. Circle packings demand to be seen; that has led to
packing algorithms, followed by experiments, then new --- often
surprising --- observations, augmented theory, more computations, on
and on. The work here was motivated by computational challenges, and
the images behind our work are produced with the open software 
package \CP, \cite{kS92}.

Step after step in this story one can observe the remarkable
faithfulness of the discrete theory to its continuous precedents so
that today one can claim a fairly comprehensive discrete world
parallel to the classical world of analytic functions (and invariably
converging to it in the limit as the combinatorics are refined).  Yet
this discrete world can never be fully comprehensive, one always faces
``discretization issues''. This paper is a preliminary description of
new machinery for addressing the principal remaining gap in the
foundation of discrete function theory, the existence and uniqueness
of discrete meromorphic functions.  The sphere is a difficult setting
for circle packing. On the practical side, there is no known algorithm
for computing circle packings {\sl in situ}, restricting the
experimental approach; essentially all circle packings on $\bP$ have
been obtained {\sl via} the stereographic projection of hyperbolic or
euclidean packings.  More crucially, the compactness of the sphere
brings conformal rigidity, with topologically mandated branching and
no boundary to provide maneuvering room.

Branching difficulties are the discretization issue we address here.
We introduce generalized branching, which began with the thesis of the
first author, \cite{jA12}. We believe general branching will provide
the flexibility necessary to construct the full spectrum of discrete
branched mappings while keeping two main objectives at the fore: (1)
discrete analytic functions should display qualitative behaviors
parallel to their classical counterparts, and (2) discrete analytic
functions should converge under refinement to their classical
counterparts.

\section{Classical Models}\label{S:Classical}
We use three types of classical functions to motivate this work:
finite Blaschke products on the unit disc, Ahlfors functions on
annuli, and Weierstrass functions on tori. We review these in
preparation for their discrete versions.

\vspace{10pt}
\noindent{\sl Blaschke Products:} A classical finite Blaschke product
$\fB:\bD\rightarrow\bD$ is a proper analytic self map of the unit disc
$\bD$. In particular, $\fB$ has finite valence $N\ge 1$, it maps the
unit circle $N$ times around itself, and it has $N-1$ branch points in
$\bD$, counting multiplicities --- that is $\fB'$ has $N-1$ zeros in
$\bD$. The function $\fB$ is known as an $N$-fold Blaschke
product. Topologically speaking, $\fB$ maps $\bD$ onto an $N$-sheeted
complete branched covering of $\bD$.  The images of the branch {\sl
  points} under $\fB$ are known as branch {\sl values}.

As a concrete example, let us distinguish two points $p_1\not= p_2$ in
$\bD$. It is well known that there exists a $3$-fold Blaschke product
$\fB$ with $p_1,p_2$ as simple branch points. It is convenient to
assume a standard normalization, so by post-composing with a conformal
automorphism (M\"obius transformation) of $\bD$ we may arrange further
that $\fB(0)=0$ and $\fB(i)=i$.  This is the function we will have in
mind for discretization later.

\vspace{10pt}
\noindent{\sl Ahlfors Functions:} Our next model is defined on a
proper annulus $\Omega$. By standard conformal mapping arguments, we
may take $\Omega$ to be a {\sl standard} annulus, $\Omega=\{z:r<
|z|<1/r\}$, with $0<r<\infty$.  Designating a point $z_0\in \Omega$,
one may consider the extremal problem: maximize $|\fF'(z_0)|$ over all
analytic functions $\fF:\Omega\rightarrow \bD$.  The solution $\fA(z)$
is known to exist, is unique up to multiplication by a unimodular
constant, and is referred to as an {\sl Ahlfors} function for
$\Omega$. Ahlfors functions are also characterized, however, by their
mapping properties. They are the proper analytic mappings
$\fA:\Omega\rightarrow\bD$ which extend continuously to
$\partial\Omega$ and map each component of $\partial \Omega$ 1-to-1
onto the unit circle. Any such map will be a branched double covering
of $\bD$ with two simple branch points, $p_1, p_2\in\Omega$. It is
fundamental to function theory on $\Omega$ and is analogous to the
1-fold Blaschke products on $\bD$, i.e., M\"obius transformations.
The Ahlfors function for $\Omega$ is determined uniquely by $r$ (up to
pre- and post-composition by conformal automorphisms).

To have a concrete example in mind for discretization, let us suppose
that $z_0$ is on the midline of $\Omega$, say $z_0=1$.  From
elementary symmetry considerations we deduce that $\fA(1)= \fA(-1)=0$
and that the branch points in $\Omega$ lie at $p_1=i$ and $p_2=-i$. A
normalization in the range, $\bD$, will put the branch values on the
imaginary axis, symmetric with respect to the origin.

\vspace{10pt}
\noindent{\sl Weierstrass Functions:} Our final model is the classical
Weierstrass function $\fW$. This is a meromorphic function mapping a
conformal torus $T$ to a branched double covering of the sphere. A
fundamental domain for $T$ is the parallelogram in $\bC$ with corners
$0,1,\tau, 1+\tau$, where $\tau$, a complex number in the upper half
plane, is the so-called {\sl modulus} of $T$. The function $\fW$ has
four simple branch points at $0,1/2,\tau/2,$ and $(1+\tau)/2$ and is
determined uniquely by $\tau$ (up to pre- and post-composition by
conformal automorphisms).

\vspace{10pt} Note that while all three classes of functions are
characterized by their topological mapping properties, only with the
Blaschke products do we get any choice in the branch points --- for
Ahlfors and Weierstrass functions, branch point locations are (up to
normalization) forced on us by the conformal geometry of the domain.

\section{Discrete Versions}\label{S:DiscVersions}
We will now describe and illustrate discrete versions of these
classical functions. We assume a basic familiarity with circle
packing, as presented in \cite{kS05} for example. However, a brief
overview might help, and with the images here should aid the
intuition, even for those not familiar with details.

A {\sl discrete analytic function} is a map between circle
packings. The domain, rather than being a Riemann surface, will now be
a triangulated topological surface with combinatorics encoded as a
simplicial 2-complex $K$: thus, we will be selecting $K$ to be a
combinatorial disc, a combinatorial annulus, or a combinatorial torus,
as appropriate. A {\sl circle packing} for $K$ is a configuration $P$
of circles, $P=\{c_v\}$ with a circle $c_v$ associated with each
vertex $v$ of $K$. The circle packing may live in the euclidean plane,
$\bC$, in the hyperbolic plane, represented as the unit disc $\bD$, or
on the Riemann sphere, $\bP$. The only requirements are that whenever
$\langle v,w\rangle$ is an edge of $K$, then circles $c_v,c_w$ must be
(externally) tangent, and when $\langle v,u,w\rangle$, is an oriented
face of $K$, then the circles $c_v,c_u,c_w$ must form an oriented
triple of mutually tangent circles. The {\sl carrier} of $P$, denoted
carr$(P)$, is the polyhedral surface formed by connecting the centers
of tangent circles with geodesic segments; that is, carr$(P)$ is an
immersion of the abstract triangulation $K$ as a concrete triangulated
surface.

At the foundation of the theory is the fact that each complex $K$ has
a canonical {\sl maximal packing} $P_K=\{C_v:v\in K\}$. This is a
univalent circle packing, meaning the circles have mutually disjoint
interiors, which fills $\bD$, a conformal annulus, or a conformal
torus, as the case may be. The packing $P_K$ serves as the domain for
discrete analytic functions associated with $K$. The image will be a
second circle packing $P$ for $K$ which lies in $\bD$ for discrete
Blaschke products and discrete Ahlfors functions, or on the sphere
$\bP$ for discrete Weierstrass functions. The discrete analytic
function, then, will be the map $\ff:P_K\rightarrow P$ which
identifies corresponding circles. (One may also treat $\ff$ as a
topological mapping $\ff:\text{carr}(P_K)\rightarrow \text{carr}(P)$
by mapping circle centers to circle centers and extending {\sl via}
barycentric coordinates to edges and faces.)

\vspace{10pt} We are now ready for the discrete constructions. Central
to our work is the issue of branching, as we will see in this first
discrete example.

\subsection{Discrete Blaschke Product}\label{SS:DiscBl}
In a sense, discrete function theory began with the introduction of
discrete Blaschke products; see \cite{tD95} and \cite[\S13.3]{kS05}.
The construction here will serve to remind the reader of basic
notation and terminology while providing an example directly pertinent
to our work.

A discrete finite Blaschke product $\fb$ is illustrated in \F{discBl},
with the domain circle packing $P_K$ on the left and the image circle
packing $P$ on the right, both in $\bD$. There is nothing special in
the underlying complex $K$, a combinatorial disc --- it is just a
generic triangulation of a topological disc, though there are minor
combinatorial side conditions to avoid pathologies.

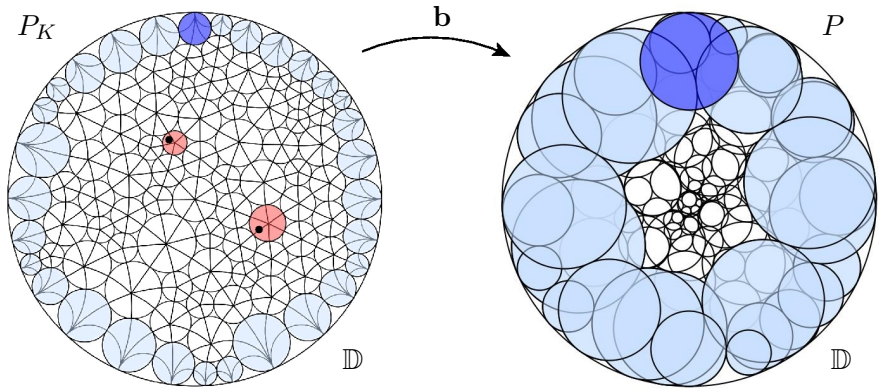
\begin{figure}[h]
\begin{overpic}[width=.95\textwidth
%,grid,tics=10
]{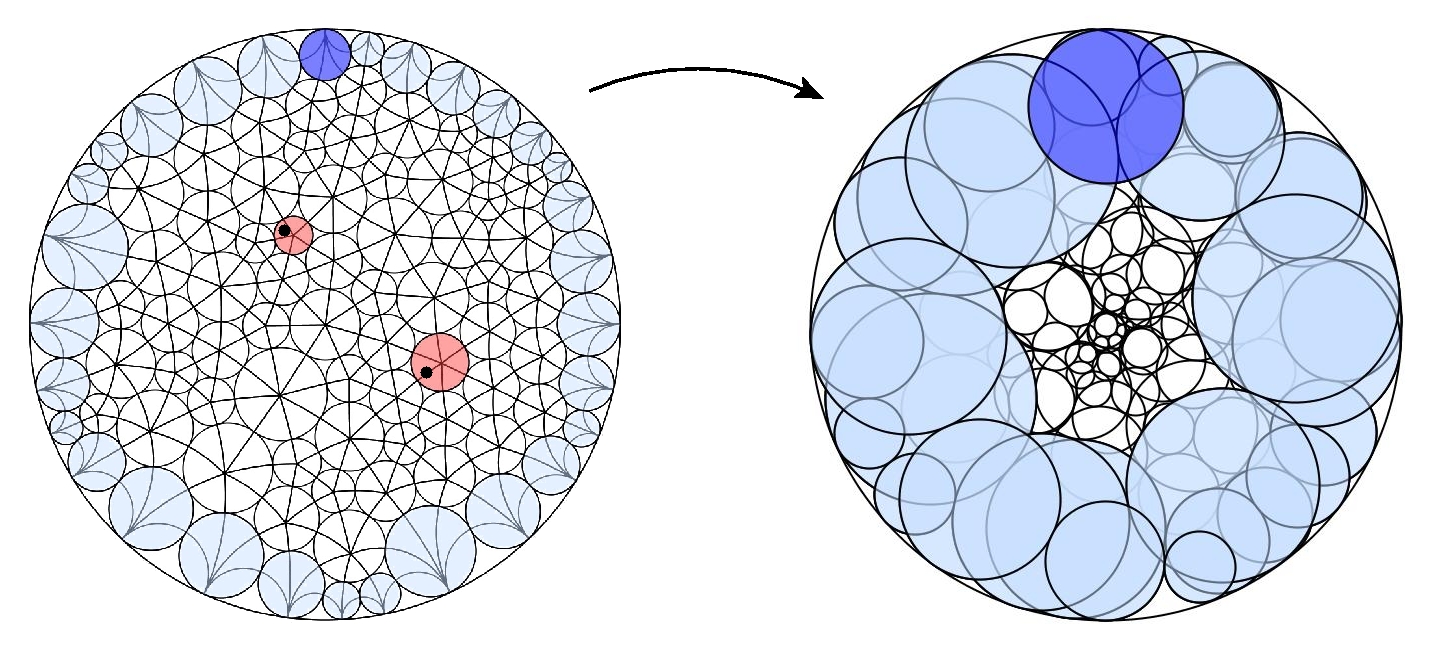}
\put (49,42) {$\fb$}
\put (39,4) {$\bD$}
\put (93,4) {$\bD$}
\put (3,41) {$P_K$}
\put (92,41) {$P$}
\end{overpic}
  \caption{A $3$-fold discrete Blaschke
    product $\fb$, domain and range.}
  \label{F:discBl}
\end{figure}

Begin with the domain packing for $\fb$ on the left, the maximal
packing $P_K=\{C_v:v\in K\}$. The boundary circles are horocyles
(euclidean circles internally tangent to $\partial\bD$). A designated
interior vertex $\alpha$ has its circle $C_{\alpha}$ centered at the
origin and a designated boundary vertex $\gamma$ has its circle
$C_{\gamma}$ centered at $z=i$; the latter appears here as dark
blue. The classical Blaschke product $\fB$ discussed earlier involved
branch points $p_1,p_2$; we assume these are the two black dots in the
domain. To mimic this, we have identified interior circles
$C_{v_1},C_{v_2}$, red circles, whose centers are nearest to
$p_1,p_2$, respectively.

Note that the unit disc is treated as the Poincar\'e model of the
hyperbolic plane, so circle centers and radii are hyperbolic and the
carrier faces are hyperbolic triangles. The boundary circles, as
horocyles, are of infinite hyperbolic radius and have hyperbolic
(ideal) centers at their points of tangency with the unit circle. The
set of hyperbolic radii is denoted by $R_K=\{R_K(v)\}$. The existence
of $P_K$ follows from the fundamental Koebe-Andreev-Thurston Theorem,
\cite[Chp~6]{kS05}, as does its essential uniqueness up to conformal
automorphisms of $\bD$. In practice, however, it is computed based on
angle sum conditions. The {\sl angle sum} $\theta_{R_K}(v)$ at a
vertex $v$ is the sum of angles at $v$ in all the faces to which it
belongs and is easily computed from the radii $R_K$ using basic
hyperbolic trigonometry. Clearly, one must have $\theta_{R_K}(v)=2\pi$
for very interior $v$. This, along with the condition that
$R_K(w)=\infty$ for boundary vertices $w$, is enough to solve for
$R_K$.

Let us now move to the more visually challenging range packing in
\F{discBl}, denoted $P=\{c_v:v\in K\}$. This, too, is a hyperbolic
circle packing for $K$, though it is clearly not univalent. We have
arranged that the circle $c_{\alpha}$ is centered at the origin and
that the circle $c_{\gamma}$ is a horocycle centered at $z=i$, just as
in $P_K$. The boundary circles are again horocycles, and if one starts
at $c_{\gamma}$ and follows the counterclockwise chain of successively
tangent horocycles, one finds that they wraps three times about the
unit circle. This mimics the behavior of our $3$-fold classical
Blaschke product $\fB$.

The image of $P$ is a bit too fussy to show its carrier, but it is in
fact a $3$-sheeted branched surface. Hidden among the interior circles
of $P$ are the two associated with vertices $v_1,v_2$, the branch
vertices.  These circles, red in both domain and range, are difficult
to pick out, but since branching is the central topic of the paper, we
have blown up the local images at $v_1$ in \F{detailBl}. We now
describe what you are seeing.

\begin{figure}[h]
\begin{overpic}[width=.7\textwidth
%,grid,tics=10
]{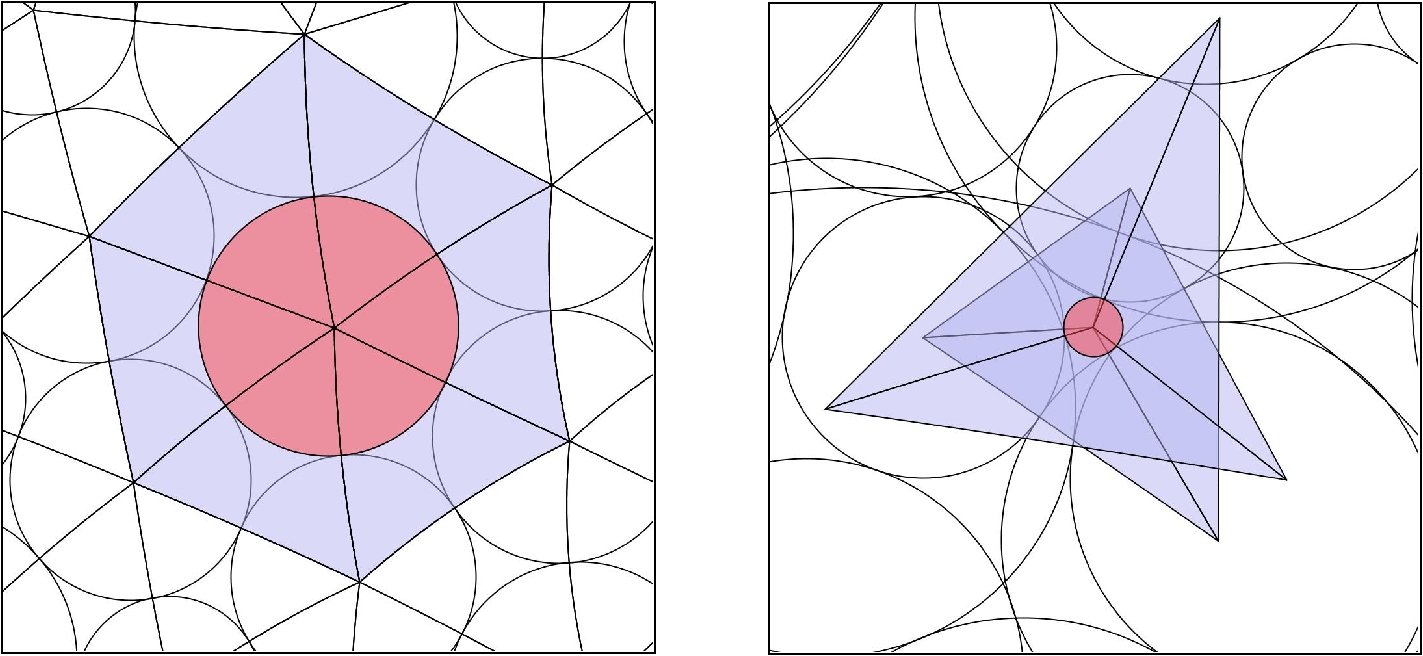}
\put (5,40) {$P_K$}
\put (92,40) {$P$}
\end{overpic}
  \caption{Isolated flowers for the branch vertex $v_1$ in
    the domain $P_K$ and range $P$ of the discrete Blaschke product
    $\fb$.}
  \label{F:detailBl}
\end{figure}

This branching will be termed {\sl traditional}; conceptually and
computationally very simple, this method has, until now, provided all
the branching for discrete function theory,
\cite[\S11.3]{kS05}. The flower for vertex $v$, the central circle (red)
and its neighboring circles (its {\sl petals}), are shown for $P_K$ on
the left and for $P$ on the right. Whereas the six petals wrap once
about $C_{v_1}$ in the domain, a careful check will show that they
wrap twice around $c_{v_1}$ in the range. If $R$ denotes the set of
hyperbolic radii for $P$, we may compute the angle sum $\theta_R(v_1)$
at $c_{v_1}$.  Expressed in terms of angle sums, the branching is
reflected in the fact that $\theta_{R_K}(v_1)=2\pi$ in the domain,
while $\theta_{R}(v_1)=4\pi$ in the range.  Mapping the faces about
$C_{v_1}$ onto the corresponding faces about $c_{v_1}$ realizes a
2-fold branched cover in a neighborhood of the center of $c_{v_1}$ ---
meaning a branched covering surface in the standard topological
sense. Similar behavior could be observed locally at the other branch
vertex, $v_2$, while at all other interior vertices the map between
faces is locally univalent.

In summary, the circle packing map $\fb:P_K\rightarrow P$ is called a
{\sl discrete finite Blaschke product} because it displays the salient
mapping features of the classical Blaschke product $\fB$: namely,
$\fb$ is a self-map of $\bD$, a 3-fold branched covering, it maps the
unit circle $3$ times about itself, and it harbors two interior branch
points.  We have even imposed the same normalization, $\fb(0)=0$ and
$\fb(i)=i$.  Additional features of such discrete analytic functions
are developed in the relevant literature: Note in \F{detailBl} how
much the circles for a branch vertex shrink under $\fb$; this ratio of
radii mimics the vanishing of the derivative at a branch point. Note
in \F{discBl} how $\fb$ draws the interior circles together; this is
the discrete hyperbolic contraction principle. Note that the circles
for $\alpha$ are centered at the origin in both $P_K$ and $P$, but the
latter is much smaller: this reflects the discrete Schwarz Lemma. On
the other hand, the horocycle associated with $\gamma$ (blue) is much
larger in $P$ than in $P_K$, reflecting the behavior of angular
derivatives at the boundary. Discrete analytic function theory is rife
with such parallel phenomena for a wide variety of situations,
including the Ahlfors and Weierstrass examples to come.

\subsection{Discrete Ahlfors Function}\label{SS:DiscAhlfors}
We build a clean example that mimics the classical Ahlfors function
$\fA$ described earlier.  Our complex $K$ triangulates a topological
annulus. Its maximal packing $P_K$ is represented in \F{AnnulusK}.

\begin{figure}[h]
\begin{overpic}[width=.6\textwidth
%,grid,tics=10
]{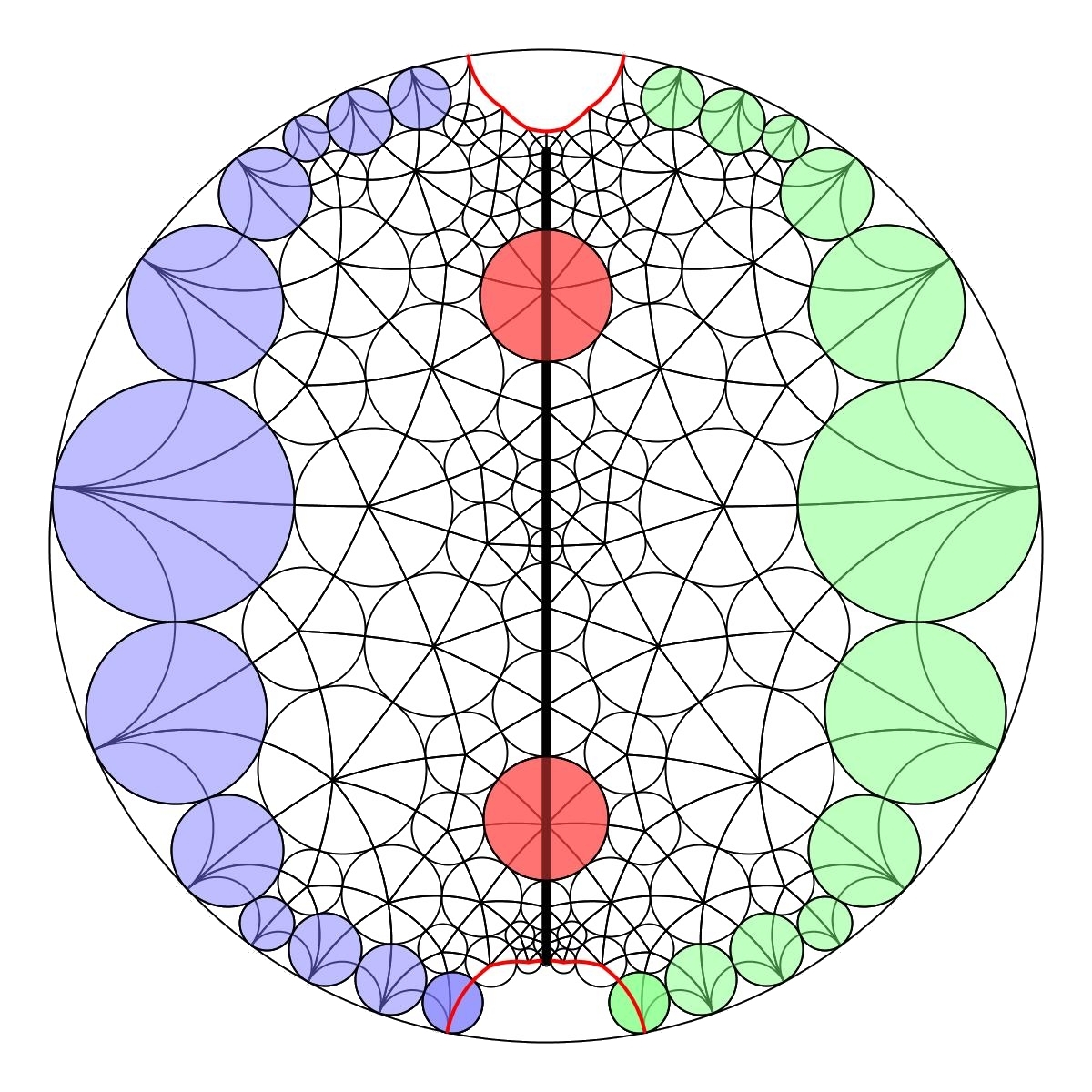}
\put (87,8) {$\bD$}
\put (4,86) {$P_K$}
\end{overpic}
  \caption{Maximal packing $P_K$ for a combinatorial annulus
    $K$, represented in a fundamental domain in $\bD$.}
  \label{F:AnnulusK}
\end{figure}

A bit of explanation may help here: The maximal packing actually lives
on a conformal annulus $\mathcal A$, with circles measured in its
intrinsic hyperbolic metric.  However, as $\bD$ is the universal cover
of $\mathcal A$, we can lift the packing to lie in a fundamental
domain within $\bD$ --- that is what we see in \F{AnnulusK}. The
boundary edges in red represent the lifts of a cross-cut of $\mathcal
A$ and are identified by the hyperbolic M\"obius transformation $\gamma$
of $\bD$ which generates the covering group for $\mathcal A$. Applying
$\gamma$ to the circles of \F{AnnulusK}, one would get new circles
which blend seamlessly along the cross-cut.

We have chosen $K$ with foresight, as it displays two particularly
helpful symmetries. The line in the center of \F{AnnulusK} marks the
combinatorial midline of the annulus: $K$ is symmetric under
reflection in this.  Moreover, there is an order two translational
symmetry along this midline; the automorphism $\sqrt{\gamma}$ will
(modulo $\gamma$) carry $P_K$ to itself. Topology demands, as with the
classical Ahlfors function $\fA(z)$, that we have two simple branch
points. Choose the midline vertices $v_1$ and $v_2$, their circles are
red in \F{AnnulusK}; these two are fixed by the reflective symmetry
and interchanged by the translational symmetry. Prescribing
traditional branching at $v_1,v_2$ results in the branched circle
packing, $P$, of \F{AnnulusP}. The mapping $\fa:P_K\rightarrow P$ is
thus a discrete analytic function from $\mathcal A$ to $\bD$.

\begin{figure}[h]
\begin{overpic}[width=.6\textwidth
%,grid,tics=10
]{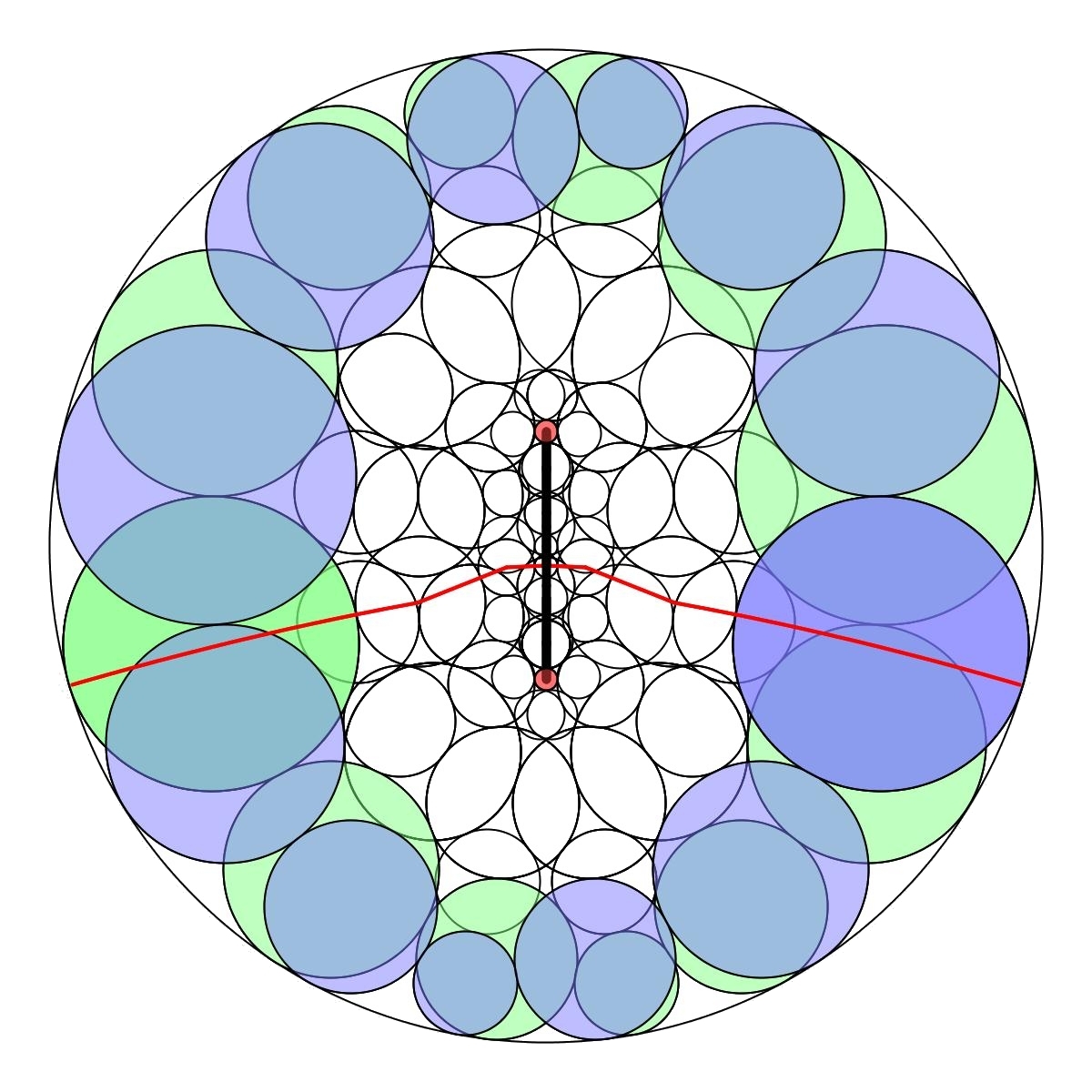}
\put (85,10) {$\bD$}
\end{overpic}
  \caption{The branched packing $P$ for combinatorial annulus $K$.}
  \label{F:AnnulusP}
\end{figure}

Due to its mapping properties, we refer to $\fa$ as a {\sl discrete
  Ahlfors function}. In particular: The boundary circles of $P_K$ are
horocycles; in \F{AnnulusK}, those on one boundary component are blue,
those of the other, green. We would expect the boundary circles of $P$
to be horocycles as well, meaning that $\fa$ maps each boundary
component to the unit circle. With a careful look in \F{AnnulusP}, one
can disentangle the closed chain of blue horocycles reaching once
around the unit circle and the second closed chain of green horocycles
doing the same. The branch circles, $C_{v_1}, C_{v_2}$ in $P_K$, and
their images, $c_{v_1},c_{v_2}$ in $P$, are red. We have normalized by
applying an automorphism to $\bD$ that centers $c_{v_1}$ and $c_{v_2}$
on the imaginary axis and symmetric with respect to the origin. Thus,
$P$ represents in a discrete way a double covering of $\bD$ branched
over two points. These are all hallmarks of the image of an Ahlfors
function and mimic the classical function $\fA$. For reference, in $P$
we have drawn in red the edges of $P$ corresponding to red cross-cut
in $P_K$.

The computation of $P$ deserves special attention. Standard Perron
methods allow one to compute a hyperbolic packing label $R$ for $K$ so
that $R(w)=\infty$ for each boundary vertex $w\in K$ and angle sums
$\theta_{R_K}(v_j)=4\pi$ for $v_1,v_2$.  There is nothing special in
computing $R$. There is a second step, however: with $R$ in hand, one
then lays out the circles in sequence and normalizes to get the
packing $P$ of \F{AnnulusP}. But {\sl why} does this second step work
so nicely?  In circle packing, the laying out of circles is akin to
analytic continuation of an analytic function element, and since $K$
is an annulus, its fundamental group is generated by some simple,
closed, non\-null\-homotopic loop $\Gamma$. Analytic continuation
along $\Gamma$ would generically lead to a non-trivial holonomy: that
is, given a function element $\mathfrak f$ defined at a point of
$\Gamma$, one would anticipate a non-trivial automorphism $m$ of $\bD$
so that analytic continuation of $\mathfrak f$ about $\Gamma$ would
lead to a new element $m(\mathfrak f)$, $m(\mathfrak f)\not= \mathfrak
f$.  In discrete terms, after laying out the circle $c_v$ for some
vertex $v$ of $\Gamma$, and then laying out successively tangent
circles for the vertices along $\Gamma$, one would not expect that
upon returning to $v$ one would lay out the same circle
$c_v$. Generically, there is a non-trivial automorphism $m$ so that
upon returning to $v$ one lays out $m(c_v)\not= c_v$. As it happens
here, things work out because of the symmetries built into $K$ --- the
holonomy $m$ is trivial, so the layout process results in a coherent
branched circle packing $P$.  The holonomy issue is key to later
considerations.

\subsection{Discrete Weierstrass Function}\label{SS:DiscWeier}
For this example our complex $K$ triangulates a topological torus.
Its maximal packing $P_K$ is shown in \F{WeierP}. Here again the
maximal packing actually lives in a conformal torus $\mathcal T$ with
its intrinsic euclidean metric. As $\bC$ is the universal cover of
$\mathcal T$, we may lift the packing to $\bC$, and this lift is what
we see on the left in \F{WeierP}. This packs a fundamental domain,
delineated by the color-coded edges, which represent the
side-pairings.

We again have chosen $K$ with important symmetries. The four colored
circles, symmetrically placed, are our chosen branch circles. The
image packing $P$ on the sphere is shown on the right in \F{WeierP};
two branched circles are visible on the front, the other two are (due
to a normalization) antipodal to these. The result is a discrete
analytic function $\fw:P_K\rightarrow P$ which maps $\mathcal T$ to
$\bP$, that is, a {\sl discrete meromorphic} function. Reprising its
mapping properties, we are justified in calling $\fw$ a {\sl discrete
  Weierstrass function}.

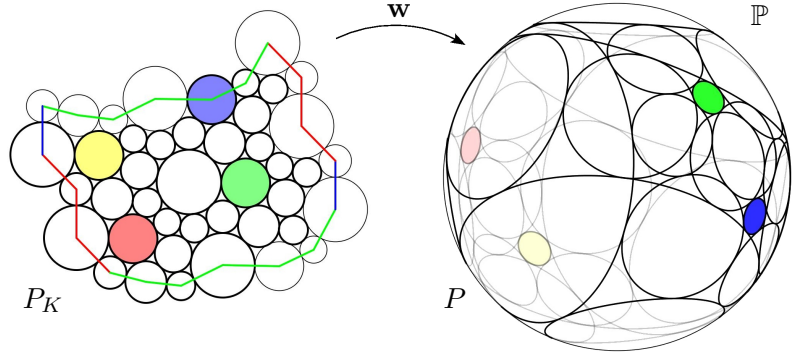
\begin{figure}[h]
\begin{overpic}[width=.85\textwidth
%,grid,tics=10
]{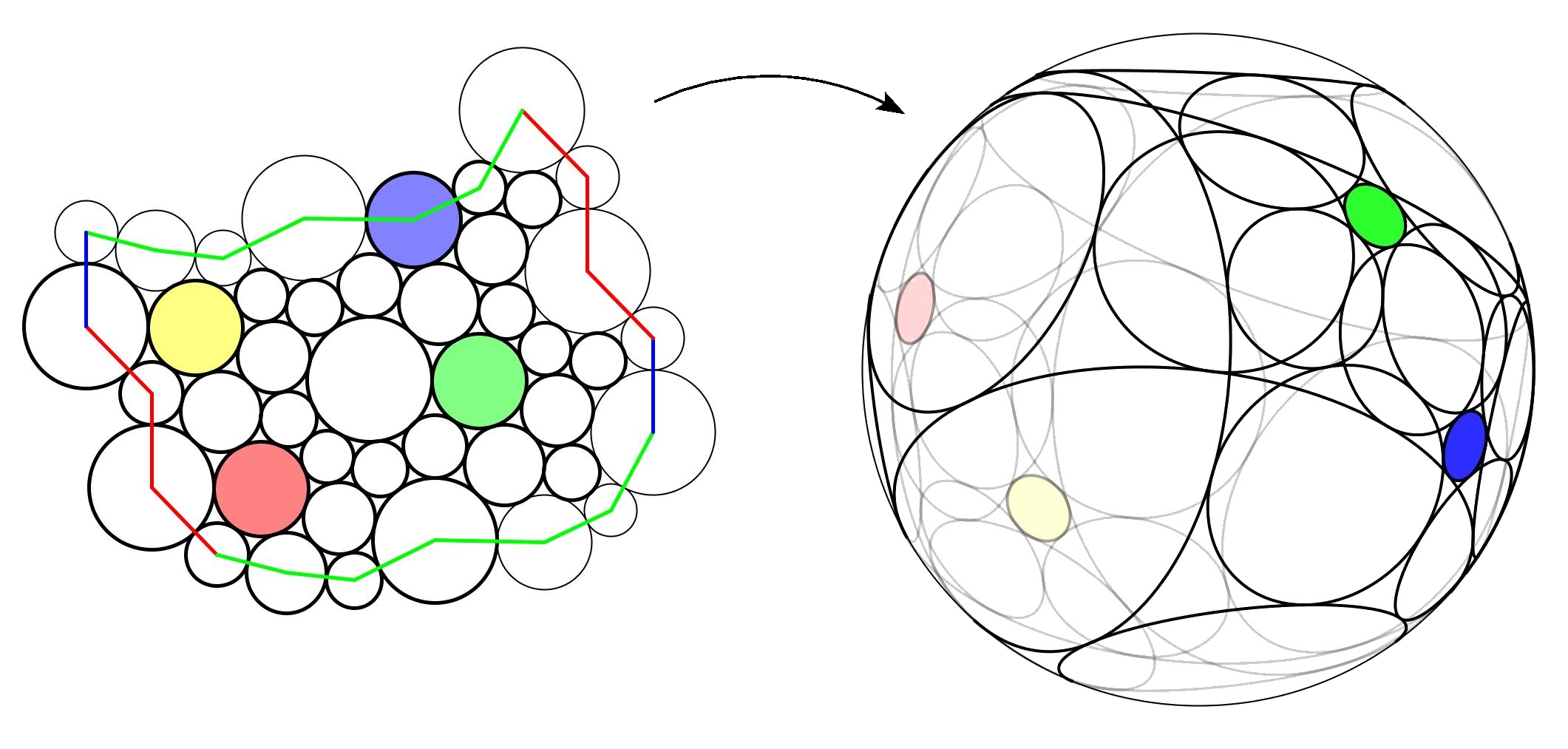}
\put (48,44) {$\fw$}
\put (3,8) {$P_K$}
\put (55,8) {$P$}
\put (93,43) {$\bP$}
\end{overpic}
  \caption{$P_K$ packs the fundamental domain for a
    combinatorial torus $K$ and shows four designated branch vertices
    and color-coded side pairings. $P$ is the branched image packing
    for $K$ on $\bP$.}
  \label{F:WeierP}
\end{figure}

There is not yet a practical circle packing algorithm in spherical
geometry, so the computation of $P$ takes a circuitous route. We
puncture $K$ at one of the intended branch vertices, say $v_4$, and
consider $K'=K\backslash\{v_4\}$. This has a single boundary
component, and the usual Perron arguments yield a hyperbolic packing
label $R'$ so that $R'(w)=\infty$ for every boundary vertex $w$ and so
that $\theta_{R'}(v_j)=4\pi, j=1,2,3$.  Since $K'$ has genus 1, its
fundamental group is again an obstruction to the layout process and a
risk for non-trivial holonomy. However, symmetry saves us once more,
and we obtain a coherent branched packing $P'$ in $\bD$ for
$K'$. Note, in particular, that the boundary circles of $P'$ are
horocycles, and topological counting arguments show that the chain of
boundary horocycles must wrap twice around the unit
circle. Stereographically projected to $\bP$, $P'$ lies in one
hemisphere. The other hemisphere, treated as the inside of a circle,
is tangent to the circles for all the former neighbors of $v_4$, so we
simply declare this to be the circle for $v_4$. The neighbors wrap
twice around, so this is the fourth branch point and, after a
normalization, we arrive at $P$. 

Note: Our methods clearly yield a
coherent branched packing $P$ in this case, and have done the same in
literally scores of similarly structured complexes.  The key seems to
lie with the symmetries in $K$ and in the branch set.  We leave this
as a {\bf Conjecture:} {\sl If $K$ is a combinatorial torus with two
  commuting translational symmetries of order two and
  $\omega=\{v_1,v_2,v_3,v_4\}$ is an orbit of vertices under these
  symmetries, then traditional branching at the points of $\omega$
  leads to trivial holonomy.}

\vspace{15pt} The good news is that we have successfully created
discrete analogues for our three classical models: Blaschke products,
Ahlfors functions, and Weierstrass functions. Let us now look into the
bad news.

\section{The Discretization Issue}\label{S:DiscIssue}
Whenever a continuous theory is discretized, whether in geometry,
topology, differential equations, or $p$-adic analysis, problems will
crop up. Replacing a continuous surface by a triangulated one, for
example, leads to combinatorial restrictions. Thus a branch vertex
must be interior and have at least 5 neighbor. We expect
this. However, we are after a starker discretization effect.  {\sl
  There are only finitely many possible locations for discrete
  branching}. Our discrete Blaschke product could not branch precisely
at the points $p_1,p_2$ prescribed for its classical model $\fB$, and
we instead chose to branch using the nearby circles $C_{v_1}$ and
$C_{v_2}$. This effect is admittedly minor --- the qualitative
behavior of the discrete function is little affected by the misplaced
branching. For the Ahlfors and Weierstrass cases, however, this
problem is existential --- discrete versions may fail to exist. We
will illustrate the problem in the Ahlfors cases --- and return to fix
it in \S\ref{S:FixAhlfors}.

\vspace{10pt} Nearly any break in the combinatorial symmetries of the
complex $K$ behind \F{AnnulusK} will cause the subsequent Ahlfors
construction to fail. Most such failures will be difficult to fix, so
we choose carefully: we make two small changes {\it via} edge flips so
that we preserve the reflective symmetry but break the translational
symmetry. The new complex will be denoted $K'$. Repeating the Ahlfors
construction from \S\ref{SS:DiscAhlfors} with $K'$ and using the same
$v_1,v_2$ as branch vertices gives the result of \F{FailedAhlfors}.

\begin{figure}[h]
\begin{overpic}[width=.58\textwidth
%,grid,tics=10
]{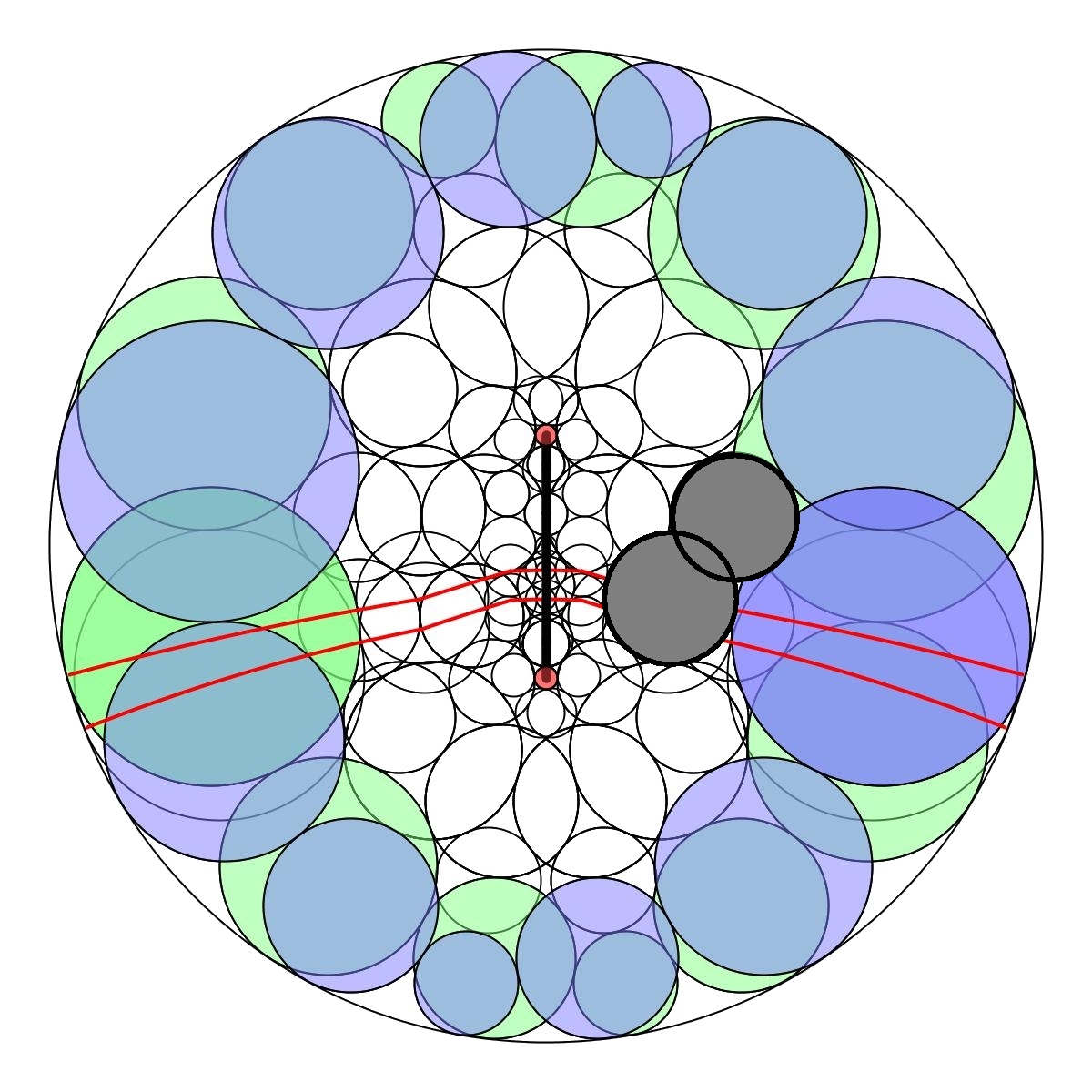}
\put (87,8) {$\bD$}
\end{overpic}
  \caption{A failed attempt at an Ahlfors function using
    traditional branching. The non-trivial holonomy shows up in
    misalignment of the cross-cuts, and the failure of the gray
    circles to be tangent to one another.}
  \label{F:FailedAhlfors}
\end{figure}

There is no difficulty in computing the branched packing label $R'$
for $K'$, however, the layout process does not give a coherent circle
packing. The problem might be difficult to see in \F{FailedAhlfors},
but look to the red edge paths, which correspond to layouts of the
cross-cut: they are no longer coincident, as they were in
\F{AnnulusP}.  One is a shifted copy of the other, reflecting a
non-trivial holonomy associated with the generator of the fundamental
group for $K'$. More precisely, there is a non-trivial hyperbolic
M\"obius transformation $m$ of $\bD$, which maps one of these red
cross-cut curves onto the other. One would have to follow things very
closely in the image to confirm the problem, but we illustrate with
the two gray circles, which are supposed to be tangent to one another.

As it happens, no matter what pair of vertices of $K'$ are chosen as
branch points, the Ahlfors construction will fail --- there will be no
coherent image packing. It has been a long road to get to this
point, but this is where our work begins: Our goal is to introduce
{\sl generalized} branching with the flexibility to make the discrete
theory whole. We will illustrate it in action in \S\ref{S:FixAhlfors}
by creating an Ahlfors function for this modified complex $K'$.

\section{Generalized Branching}\label{S:GenBranching}
Branching is perhaps most familiar in the analytic setting.  Let
$\ff:G\rightarrow \bC$ be a non-constant analytic function on an open
domain $G\subset \bC$.  Suppose $z\in G$ and $w=\ff(z)$. For
$\delta>0$, consider the disc
$D=D(w,\delta)=\{\zeta:|\zeta-w|<\delta\}$ and the component of the
preimage $U={\ff}^{-1}(D(w,\delta))$ containing $z$.  For $\delta$
sufficiently small, $U$ will be a topological disc in $G$ and the
restriction of $\ff$ to the punctured disc $U'=U\backslash\{z\}$ will
be a locally 1-to-1 proper mapping onto the punctured disc
$D'=D\backslash\{w\}$. In particular, one can prove the existence of
some $N\ge 1$ so that every point of $D'$ has $N$ preimages in
$U'$. In this analytic case, if $N>1$, then ${\ff}^{k}(z)=0$ for
$k=1,2,\cdots, N-1$ and we say that $z$ is a {\sl branch point} of
order $N-1$ for $\ff$. We refer to $w$ as its {\sl branch value}.

This is, in fact, a topological phenomenon having little to do with 
analyticity: by Sto\"ilow's Theorem, \cite{sS38}, the same local behavior 
occurs whenever the map $\ff$ is an open, continuous, light-interior 
mapping. In particular, this applies to our maps between the carriers 
of circle packings. One sees it on display for the traditional branch 
point illustrated in \F{detailBl}.

We can set the stage for generalized branching by simply enlarging the
singleton set $\{z\}$ for a branch point to a compact topological disc
$H$. If $H$ is small enough, then the mapping behavior in the
neighborhood of $H$ is unchanged: that is, $\ff$ will be a locally
1-to-1 proper map of the annulus $U'=U\backslash H$ onto the annulus
$D'=D\backslash f(H)$ of valence $N$. When $N>1$ we will say that
$\ff$ has {\sl generalized branching} of order $N-1$ in $H$. The point
is that the branching is reflected in the mapping behavior between the
annuli $U'$ and $D'$, even if the precise location of that branching
is hidden within the hole $H$.

Let us apply this notion to the classical Blaschke product $\fB$
discussed earlier. About each of its branch points $p_j,j=1,2,$ we can
choose a small compact topological disc $H_j$ and open neighborhood
$U_j$ of $H_j$ so that $\fB$ has generalized branching of order $1$ in
$H_j$. Making $H_1,H_2$ smaller, if necessary, we may assume $U_1$ and
$U_2$ have disjoint closures. This leaves a triply connected open set
$\Omega=\bD\backslash\{\overline{U}_1\cup \overline{U}_2\}$. The
restriction of $\fB$ to $\Omega$ is a $3$-valent map onto $\bD$ which
maps $\partial U_j$ to a some curve about the branch value
$\fB(p_j),j=1,2$. If we were to perturb $p_j$ to a new point $\tilde
p_j$ within $H_j, j=1,2$, then the associated finite Blaschke product
$\widetilde{\fB}$ would be qualitatively indistinguishable from $\fB$
on $\Omega$: that is, it is very difficult to discern where the branch
points actually lie in $H_1,H_2$. This is the type of flexibility we
need for discrete finite Blaschke products and motivates our strategy.

Given $K$ and its max packing $P_K$, we choose interior vertices
$v_1,v_2$ so their circle centers are near $p_1,p_2$. Choose small
combinatorial neighborhoods $H_1,H_2$ of $v_1,v_2$ and define
$L=K\backslash\{H_1\cup H_2\}$, analogous to the open set $\Omega$
earlier. Requiring simple branching at $v_1$ and $v_2$ leads to the
discrete Blaschke product $\fb$ we discussed earlier. However, we have
developed machinery, discussed in the next section, that allows us to
modify the combinatorics and packing parameters inside $H_1$ and
$H_2$. Patching these new combinatorics into $L$ gives a new
combinatorial disc $\wK$ on which we define a new discrete
Blaschke product $\tilde{\fb}$. The parameters involved allow us to
perturb the apparent branch locations.  In other words, just as with
$\fB$ and $\widetilde{\fB}$ defined on $\Omega$, both $\fb$ and
$\tilde{\fb}$ are defined on $L\subset K\cap \wK$ and are
qualitatively indistinguishable there.

Our global intention is to make adjustments in the small locales $H_1$
and $H_2$ so that $\tilde{\fb}$ behaves like a discrete Blaschke
product having branch points precisely at $p_1,p_2$. Mimicking this
individual Blaschke product $\fB$ may seem to be a lot of effort for
little gain. However, if one thinks more broadly of the family of
Blaschke products parameterized continuously by $p_1$ and $p_2$, the
goal of continuously parameterized discrete versions makes more
sense. It also makes more sense when the very existence of the
discrete versions depends on this added flexibility, as with our
broken Ahlfors example. Let us now describe the mechanics.

\section{Local Mechanics}\label{S:Local}
We describe discrete generalized branching which takes two forms,
termed {\sl singular} and {\sl shifted} branching. Each involves
identifying a {\sl black hole} $H$, a small combinatorial locale to
support the branching, and its {\sl event horizon} $\Gamma=\partial
H$, the chain of surrounding edges.  Outside of the event horizon, our
circle packing mappings are defined in the usual way, so that in an
annulus about the black hole one may observe the typical topological
behavior described earlier. Adjustments hidden inside the black hole,
however, allow our mapping to simulate simple branching at various
points.

\subsection{Background}\label{SS:Background}
We have recalled some circle packing mechanics, but as our work
involves new features, we review the basics.

A complex $K$ is assumed to be given. The fundamental building blocks
of $K$ are its {\sl triangles} and {\sl flowers}. The triangles are
the faces $\langle v,u,w\rangle$.  The flowers are sets
$\{v;v_1,v_2,\cdots,v_{n+1}\}$ where $v$ is a vertex and
$v_1,\cdots,v_{n+1}$ is the counterclockwise list of neighbors in $K$.
These neighbors, the {\sl petals}, define the fan of faces containing
$v$. Here $n$ is the number of faces, so when $v$ is interior, then
$v_{n+1}=v_1$.

In talking about a circle packing $P$ for $K$, the radii and centers
are, of course, the ultimate target. However, proofs of existence and
uniqueness (and computations) depend on the standard Perron methods
first deployed in \cite{BSt91a}. Given $K$, the fundamental data lies
in three lists: the {\sl label} $R=\{R(v):v\in K\}$, edge overlaps
$\Phi=\{\Phi(e): e=\langle u,v\rangle\}$, and target angle sums
$A=\{A(v):v\in K\}$. Each will require some extension.

\begin{itemize}
\item{{\bf Labels:} The labels $R(v)$ are putative radii (they become actual
  radii only when a concrete packing is realized).}

\vspace{5pt}
\item{{\bf Overlap Angles:} For an edge $e=\langle v,w\rangle$ of $K$,
  the overlap $\Phi(e)$ represents the desired (external) angle
  between the circles $c_v,c_w$ in $P$. Interest is often in
  ``tangency'' packings; in this case, $\Phi$ is identically zero and
  hence does not appear explicitly. However, from Thurston's first
  introduction of circle packing, non-tangency packings were included
  and we need them here.}

\vspace{5pt}
\item{{\bf Target Angle Sums:} Given $R$ and $\Phi$, one can readily
  compute for any triangle $\langle u,v,w\rangle$ the angle which
  would be realized at $v$ if a triple of circles with the given
  labels (as radii) and edge overlaps were to be laid out.  The {\sl
    angle sum} $\theta_{R,\Phi}(v)$ is the sum of such angles for all
  faces containing $v$. The target angle sum, $A(v)$ is the
  intended value for $\theta_{R,\Phi}(v)$. It is typically prescribed
  only when $v$ is interior, and then must be an integral multiple of
  $2\pi$, $A(v)=2\pi k$; this is precisely the result when petal
  circles $c_{v_1},c_{v_2},\cdots,c_{v_n}$ wrap $\frac{A(v)}{2\pi}=k$
  times about $c_v$.}
\end{itemize}

A circle packing for $K$ is computed by finding a label $R$, termed a
{\sl packing label}, with the property that $\theta_{R,\Phi}(v)=A(v)$
for every interior vertex $v$. Typically, the values $R(w)$ for
boundary vertices $w$ are prescribed in advance. With label in hand,
one can position the circles in the pattern of $K$ to get $P$. This
positioning stage is a layout process analogous to analytic
continuation for analytic functions. Only after the layout does one
finally realize circle centers. Our work will be carried out in
hyperbolic geometry, where we use the fact that boundary radii may be
infinite when associated with horocycles. The various existence,
uniqueness, and monotonicity results needed for our applications would
hold in euclidean geometry as well.

The Perron method for computing a packing label proceeds {\sl via}
{\sl superpacking} labels, that is, labels $R$ for which the
inequality $\theta_{R,\Phi}(v)\le A(v)$ holds for all interior $v$ and
which has values no less than the designated values at boundary
vertices. It is easy to show that the family of superpacking labels is
nonempty and that the packing label is the family's infimum. This
infimum may be approximated to any desired accuracy by an iterative
adjustment process --- this is basically how \CP\ computations are
carried out.  The following condition ($\star$) is required to ensure
non-degeneracy: If $\{e_1, e_2,\cdots,e_k\}$ is a simple closed
edge path in $K$ which separates some edge-connected non-empty set $E$
of vertices from $\partial K$, then the following inequality must hold
\begin{equation}\tag{$\star$}\label{E:blackhole}
\sum_1^k(\pi-\Phi(e_j))\ge 2\pi+\sum_{v\in E}(A(v)-2\pi).
\end{equation}

Our work here requires the following extensions to the given
data:  

\begin{itemize}
\item{{\bf Zero Labels:} We will introduce situations in which labels
  for certain interior vertices go to zero, corresponding with circles
  that in the final configuration have degenerated to points, namely
  to their centers. Zero radii actually fit quite naturally into the
  trigonometric computations, but we will only encounter them for
  isolated vertices.}

\vspace{5pt}
\item{{\bf Deep Overlaps:} When introducing circle packing, Thurston
  included specified overlaps $\Phi(e)$, as we do. In general,
  however, the restriction $\Phi(e)\in [0,\pi/2]$ is required for
  existence. We will allow {\sl deep overlaps}, that is overlaps in
  $(\pi/2,\pi]$. Note that overlaps may already be specified as part
  of the original packing problem under consideration, but these will
  remain in the range $[0,\pi/2]$. It is only in the modifications
  within black holes that deep overlaps may be needed, and these will
  carry clear restrictions.}

\vspace{5pt}
\item{{\bf Branching:} Traditional branching, described earlier in the
  paper, is associated with target angle sums $A(v)=2\pi k$ for $k\ge
  2$. These are subject to the condition ($\star$) noted above, which
  concerns interactions of combinatorics and angle sum
  prescriptions. It traces to the simple observation that it takes at
  least 5 petal circles to go twice around a circle. The tight
  conditions emerged first in work on branched tangency packings in
  \cite{tD93} and \cite{pB93}. These were modified to incorporate
  overlaps in \cite{BS96}; Condition ($\star$) parallels the conditions
  there while allowing equality, which is associated with zero labels in
  black holes, as we see shortly.}

\end{itemize}

The monotonicities behind the Perron arguments depend on our ability
to realize any face $\langle u,v,w\rangle$ with a triple of circles
$\{c_v,c_u,c_w\}$ having prescribed radii and overlaps.  To include
deep overlaps and zero labels, it is relatively easy to see that some
side conditions on $\Phi$ are required. What we need is given in the
following lemma, a minor extension of the hyperbolic results in
\cite{jA12}.

\begin{lemma}\label{L:Ashe} Given three hyperbolic radii, $r_1,r_2,r_3$,
at least two of which are non-zero, and given three edge overlaps
$\phi_{12},\phi_{23},\phi_{31}\in[0,\pi]$ satisfying
$\phi_{12}+\phi_{23}+\phi_{31}\le \pi$, there exists a triple $\langle
c_1,c_2,c_3\rangle$ of circles in $\bD$ which realize the given radii
and overlaps.

The angles $\alpha,\beta,\gamma$ of the triangle $T$ formed by their
centers are continuous functions of the radii and overlaps and are
unique up to orientation and conformal automorphisms of $\bD$.
Moreover $\alpha$ is strictly decreasing in $r_1$, while area$(T)$ is
strictly increasing in $r_1$. Likewise, $\beta$ (resp. $\gamma$) is
strictly increasing in $r_1$ (assuming $r_2$ (resp. $r_3$) is finite).
\end{lemma}

In our generalized branching, zero labels and deep overlaps are
temporary devices only within black holes; we modify the combinatorics 
and set overlap parameters in there to control apparent branch
locations. The results, however, are then used to layout a circle
packing $P$ for the original complex $K$; $P$ itself does not involve
any zero labels or deep overlaps, and aside from ambiguity about one
circle in the shifted branching case, $P$ is a normal circle packing
configuration.

We conclude these preparations by noting the two conditions which are
necessary to guarantee existence and uniqueness of the
packings. Namely, condition ($\star$) and this condition ($\star\star$)
\begin{equation}\tag{$\star\star$}\label{E:DeepCondition}
\Phi(e_1)+\Phi(e_2)+\Phi(e_3)\le \pi\ \text{if edges }
\{e_1,e_2,e_3\}\ \text{form a face of }K.
\end{equation}
With this, we may now describe our two discrete generalized branching
mechanisms.

\subsection{Singular Branching}\label{SS:SingBranching}
Singular branching is used to simulate a branch point lying in an
interstice of $P_K$. The interstice is defined by a face $\langle
v_1,v_2,v_3\rangle$, corresponding to red, green, and blue circles,
respectively, in our illustrations. The black hole is the union of the
target interstice and the three interstices sharing its edges.  The
combinatorics imposed and the event horizon are illustrated in
\F{SingComb}. The complex $K$, modified inside the black hole, will be
denoted $\wK$ and serves as our complex for subsequent computations.
The circles of \F{SingComb} are a device for display only and are not
part of the final circle configuration. Indeed, before computing the
circles of the branched packing we need to prescribe target angle
sums, $A$, and edge overlaps, $\Phi$.

\begin{figure}[h]
\begin{overpic}[width=.5\textwidth
%,grid,tics=10
]{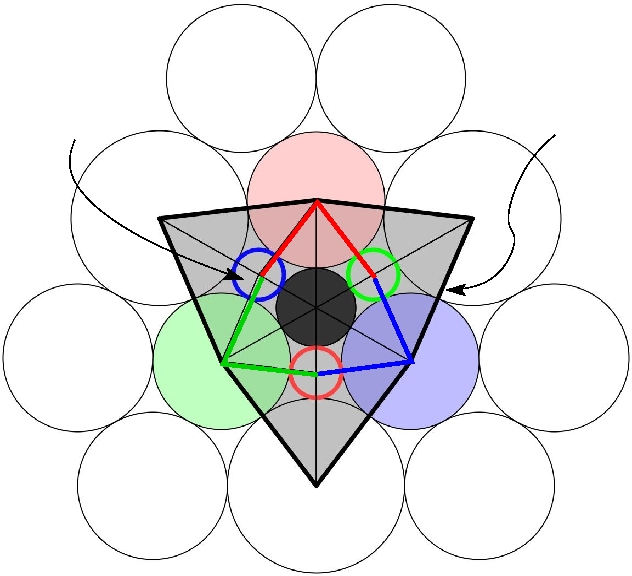}
\put (51,62) {$v_1$}
\put (29,30) {$v_2$}
\put (67,30) {$v_3$}
\put (82,77) {event}
\put (82,72) {horizon}
\put (-2,76) {chaperone}
\put (10,70) {$h_3$}
\end{overpic}
  \caption{Combinatorics for a singular black hole.}
  \label{F:SingComb}
\end{figure}
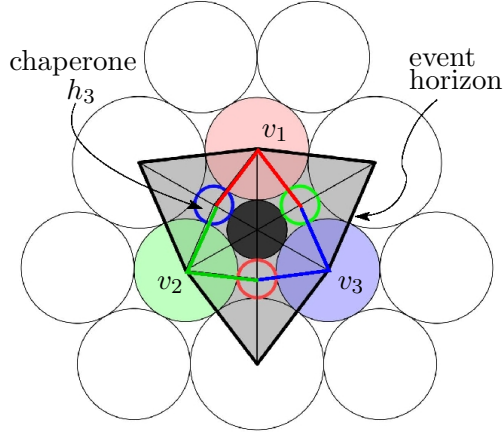

Interior to the event horizon we have introduced 4 additional
vertices. Three of these, $h_1,h_2,h_3$, are termed {\sl chaperones}
since they help guide the circles for $v_1,v_2,$ and $v_3$; we label
$h_3$ in \F{SingComb}. A fourth vertex $g$, in the center, is called
the {\sl fall guy}. Specify target angle sums $A(v)\equiv 2\pi$ for
all interior vertices $v\in \wK$ with the exception of $g$, setting
$A(g)=4\pi$.

Singular branching is controlled {\sl via} overlap parameters
associated with a partition of $\pi$,
$\gamma_1+\gamma_2+\gamma_3=\pi$.  For $i=1,2,3,$ the value
$\gamma_i>0$ represents the overlap angle prescribed in $\Phi$ for the
edges from $v_i$ to the chaparone circles on either side.  These three
pairs of edges are color coded in \F{SingComb}. We set $\Phi(e)= 0$
for all other edges of $\wK$.

Before describing how these parameters are chosen, observe that we are
assured of a circle packing $\wP$ for $\wK$ with label $\wR$, interior
angle sums $A$, and overlaps $\Phi$. In particular, if $\Gamma$
denotes the chain of 6 colored faces surrounding the fall guy, $g$,
then condition ($\star$) holds whenever the angle sum prescription $A(g)$
satisfies $A(g)\le 4\pi$, with equality when $A(g)=4\pi$. Traditional
Perron and layout arguments imply the existence and uniqueness of the
circle packing $\wP$ in which the circle for $g$ has radius zero.  An
example of the result is illustrated in \F{SingBranch}. For this we
set roughly $\gamma_1=0.22\pi,\gamma_2=0.40\pi$, and
$\gamma_3=0.37\pi$.

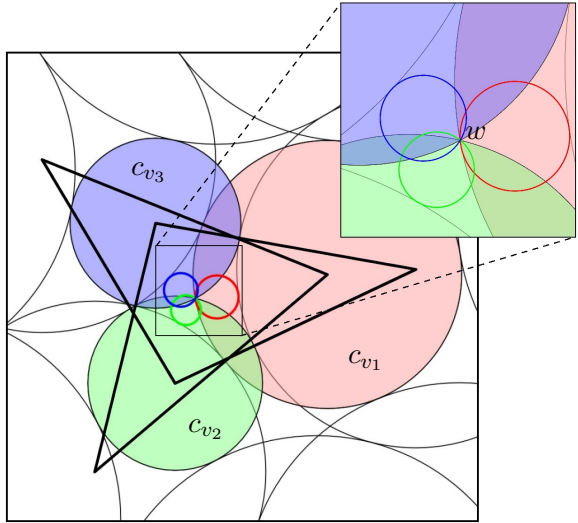
\begin{figure}[h]
\begin{overpic}[width=.6\textwidth
%,grid,tics=10
]{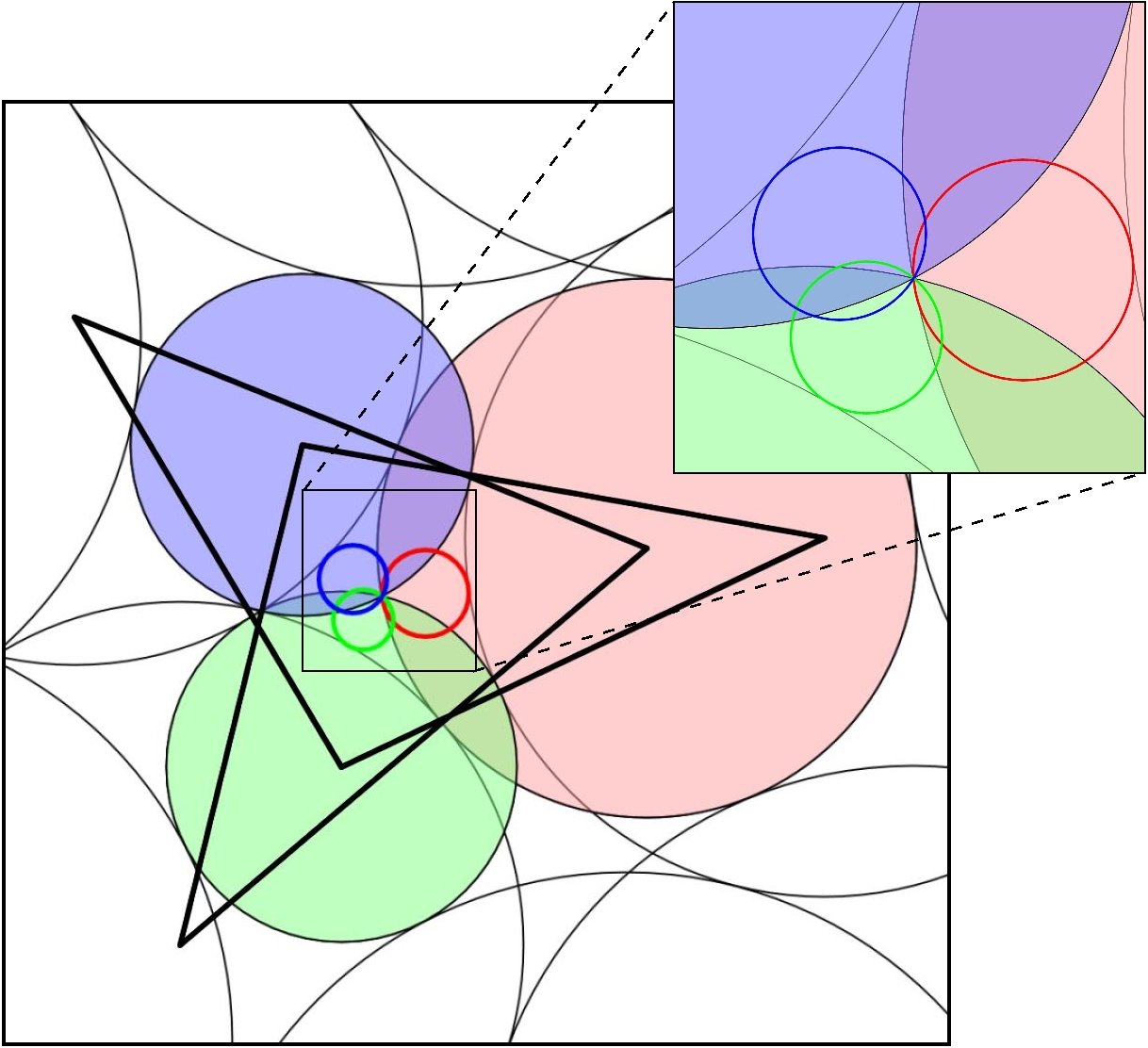}
\put (60,28) {$c_{v_1}$}
\put (32,16) {$c_{v_2}$}
\put (22,61) {$c_{v_3}$}
\put (80.5,67.25) {$w$}
\end{overpic}
  \caption{Image circle packing in the neighborhood of 
singular branching, and detail zoom.}
  \label{F:SingBranch}
\end{figure}

This image takes some time to understand. The circle for $g$ has
degenerated to a point, the branch value, which is at the common
intersection point of the circles for $v_1,v_2,v_3$ and also for
chaperones $h_1,h_2,h_3$; it is labeled $w$ in the detail zoom. The
branching is confirmed in the larger image by observing how the event
horizon wraps twice about the branch value.  If we disregard the
chaperones and the fall guy, the remaining circles of $\wP$ realize a
tangency circle packing for the original complex $K$. That is, the
black hole structure was needed only to guide the layout of the
original circles.

This portion of the layout can best be understood as living on a two
sheeted surface $S$ branched above $w$. Note, for instance, that the
overlap of the red and blue circles is only in their projections to
the plane: in actuality, the red part of the intersection is on one
sheet of $S$ and the blue is on the other. This shows in the
orientation of the red, green, blue, which in projection is the
reverse of their orientation in $P_K$.

Finally, what about choosing parameters $\gamma_1,\gamma_2,\gamma_3$
to get the desired branch point?  \F{SingScheme} illustrates our
scheme. We have isolated the interstice formed by circles for
$v_1,v_2,v_3$ in $P_K$. The dashed circle is the common orthogonal
circle through the intersection points and defines a disc $D$ which
will be treated as a model of the hyperbolic plane. Point $p$
indicates a location where one might wish to have branching
occur. Hyperbolic geodesics connecting $p$ to the three intersection
points on $\partial D$ determine angles $\alpha_1,\alpha_2,\alpha_3$,
indexed to correspond with the vertices $v_1,v_2,v_3$. We then define
$\gamma_j=\pi-\alpha_j,j=1,2,3$.  One has complete freedom to choose
$\gamma_1$ and $\gamma_2$ in this scheme, subject to conditions
$\gamma_1,\gamma_2>0$ and $\gamma_1+\gamma_2<\pi$. We will be seeing
examples for $p$ and the other three red branch points later, in
\F{SingMovie}.

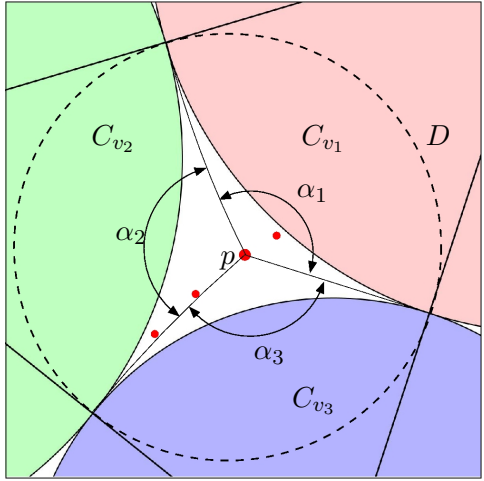
\begin{figure}[h]
\begin{overpic}[width=.5\textwidth
%,grid,tics=10
]{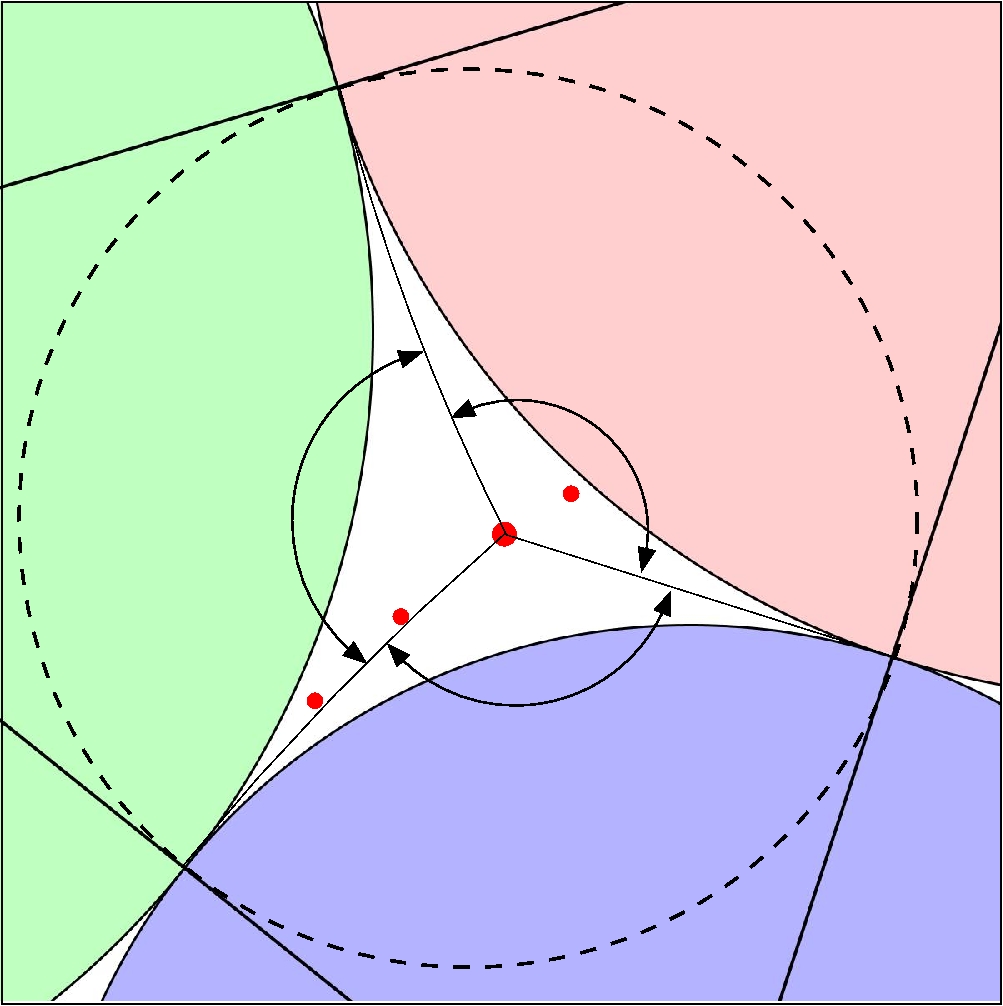}
\put (61,59) {$\alpha_1$}
\put (23,50) {$\alpha_2$}
\put (52,25) {$\alpha_3$}
\put (45,45) {$p$}
\put (88,70) {$D$}
\put (62,70) {$C_{v_1}$}
\put (18,70) {$C_{v_2}$}
\put (60,15) {$C_{v_3}$}
\end{overpic}
  \caption{The parameter scheme for singular branching.}
  \label{F:SingScheme}
\end{figure}

\subsection{Shifted Branching}\label{SS:ShiftBranching}
Shifted branching simulates a branch point lying within an interior
circle of $P_K$. Of course, when that point is the center, then
traditional branching would be the easy choice. This will be
incorporated naturally in our parameterized version, however, so we
need not separate out this case.

Suppose $v$ is the interior vertex whose circle is to contain the
shifted branch point. The black hole combinatorics shown in
\F{ShiftComb} are imposed on the flower for $v$. (Note that once
again, the circles here are used for display but are not part of our
target packing.)

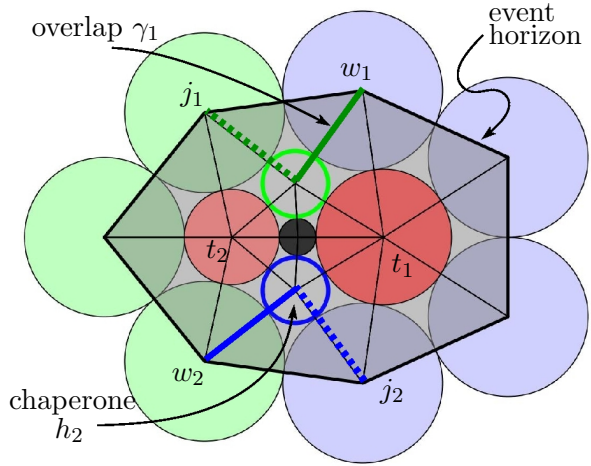
\begin{figure}[h]
\begin{overpic}[width=.6\textwidth
%,grid,tics=10
]{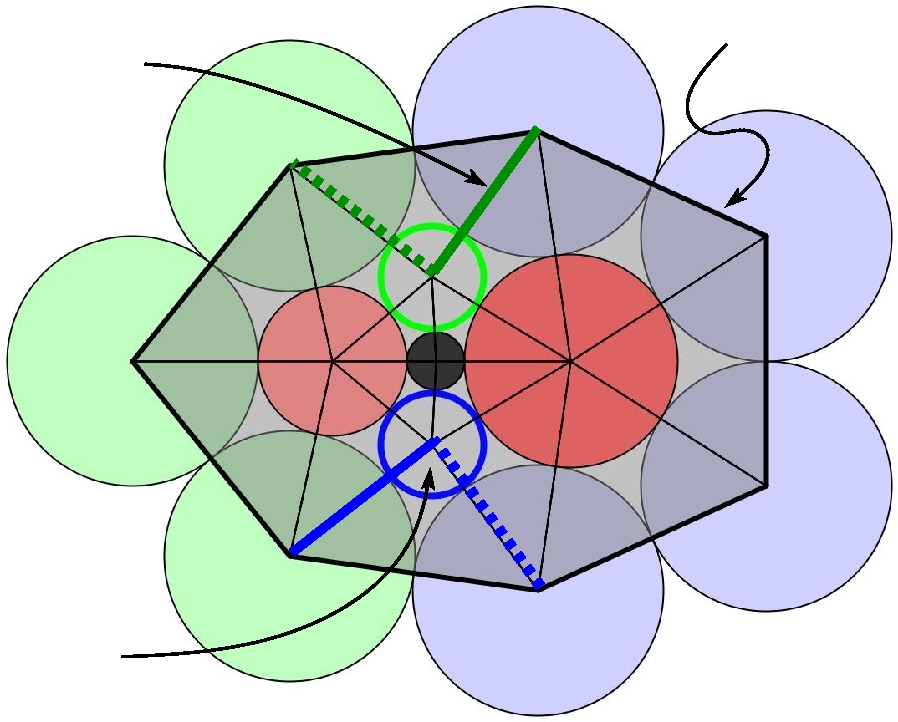}
\put (65,34) {$t_1$}
\put (32.5,36.5) {$t_2$}
\put (28,64) {$j_1$}
\put (56,68.5) {$w_1$}
\put (63,12) {$j_2$}
\put (26.5,15) {$w_2$}
\put (82,78) {event}
\put (82,74) {horizon}
\put (-2,10) {chaperone}
\put (6,5) {$h_2$}
\put (2,75.5) {overlap $\gamma_1$}
\end{overpic}
  \caption{Combinatorics for a shifted black hole.}
  \label{F:ShiftComb}
\end{figure}

The event horizon is the chain of edges through the original petals of
the flower for $v$ (seven petals, in this case, green and blue).
Interior to this horizon, we split $v$, replacing it with the twin
vertices, denoted $t_1$ and $t_2$ and corresponding to the circles in
two shades of red. We introduce two chaperone vertices $h_1,h_2$,
respectively green and blue, and a fall guy vertex $g$, black; we
label only chaperone $h_2$ in the figure.  With these combinatorics
inside the event horizon, we again have a new complex $\wK$, for which
we need target angle sums $A$ and edge overlaps $\Phi$.

Each chaperone neighbors two original petals, denoted $w_i$ and
$j_i$. The petal $j_i$ is known as the {\sl jump} circle because its
chaperone $h_j$ and an associated parameter $\gamma_i$ facilitate its
detachment from one twin and its attachment to the other.  The
parameters here are $\gamma_1$ and $\gamma_2$, chosen independently
within $[0,\pi]$, and used to define overlaps with the chaperones. In
particular, for $i=1,2,$ prescribe $\Phi(\langle
h_i,w_i\rangle)=\gamma_i$ and $\Phi(\langle h_i,j_i\rangle)=
\pi-\gamma_i$; the edges are shown as solid and dashed lines,
respectively, in \F{ShiftComb}.  The other overlaps in $\Phi$ are
zero, so Condition ($\star\star$) holds.  Target angle sums are defined as
before, namely, $A= 2\pi$ at interior vertices of $\wK$, save for the
fall guy, with $A(g)=4\pi$.

Putting aside the choice of jump circles and parameters for now, we
are assured of a circle packing $\wP$ for $\wK$ with label $\wR$,
interior angle sums $A$, and overlaps $\Phi$.  If $\Gamma$ is the
chain of edges through the four neighbors of $g$, edges for which
tangency is specified, then equality holds in condition ($\star$), so in
$\wR$ the radius of the fall guy is necessarily zero.

\F{Shift2pi} illustrates the circle packing for $\wK$ before we
prescribe the branching, in other words, with the target angle sum at
$g$ kept at $2\pi$. We abuse notation by referring to circles by their
vertex indices. The original petal circles, starting with $j_1$ and
ending at $w_2$, are shown in green: these are tangent to twin
$t_2$. Likewise, those starting at $j_2$ and ending at $w_1$ are shown
in blue: these are tangent to twin $t_1$.

\begin{figure}[h]
\begin{overpic}[width=.6\textwidth
%,grid,tics=10
]{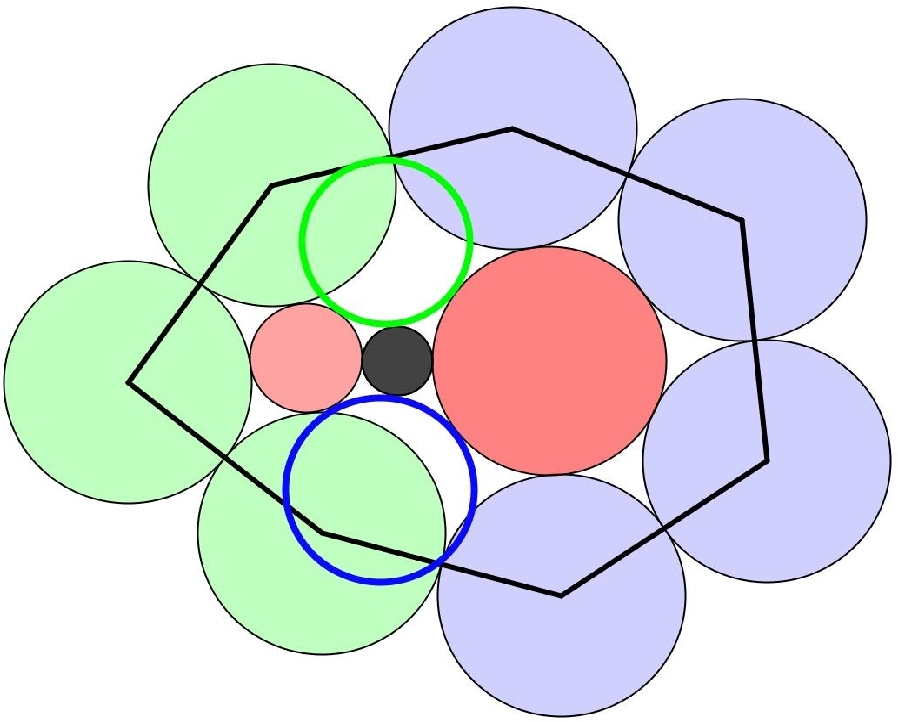}
\put (59,38) {$t_1$}
\put (33,38) {$t_2$}
\put (27,63) {$j_1$}
\put (54,68.5) {$w_1$}
\put (64,11) {$j_2$}
\put (28,17) {$w_2$}
\end{overpic}
  \caption{The jump circles and parameters are set for a
    shifted black hole before the branching is imposed.}
  \label{F:Shift2pi}
\end{figure}
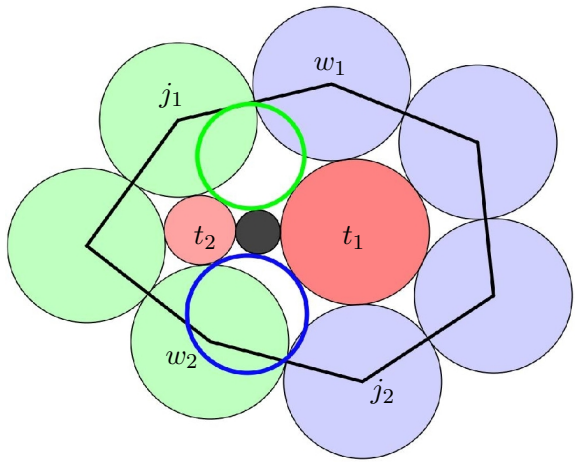

We consider the action at chaperone $h_1$. First, recall two facts:
(1) When a triple of circles has edge overlaps summing to $\pi$, then the
three share a common intersection point; and (2) when circles overlap
by $\pi$ then one is interior to the other. Here is how the machinery
works at $h_1$. The circle for $w_1$ is tangent to twin $t_1$, $j_1$
is tangent to twin $t_2$, while $h_1$ is tangent to both twins. When
$\gamma_1=0$, the overlap of $\pi$ between $h_1$ and $j_1$ forces the
jump circle $j_1$ to be tangent to $t_1$. As $\gamma_1$ increases,
however, the jump circle separates from $t_1$ until, when $\gamma_1$
reaches $\pi$, $w_1$ has been pulled in to be tangent to $t_2$. In
other words, $\gamma_1$ acts like a dial: when positive, it detaches
the jump circle from $t_1$, and as it increases, it moves the jump
further around $t_2$. The mechanism is similar for chaparone $h_2$, as
$\gamma_2$ serves to detach the jump circle $j_2$ from $t_2$ and move
it further around $t_1$.

Typical parameters $\gamma_1=0.7\pi$ and $\gamma_2=0.4\pi$ were
specified for \F{ShiftBranch}. Maintaining these while adding
branching at $g$, i.e., setting $A(g)=4\pi$, gives the configuration
of \F{ShiftBranch}.  As usual with branching, the image is rather
difficult to interpret, so we point out the key features: The twin
circles and chaperones are all tangent to $g$, and the radius for $g$
is zero, so these four circles meet at a single point.  The twin
circles (red) are nested, as are the chaperones (green and blue). The
branch value is the white dot in the detail zoom, at the center of the 
small twin circle and labeled $w$; we explain this shortly. To confirm the
topological behavior of generalized branching, note that the circles
for the original petals of $v$ wrap twice around $w$ --- just follow
the image of the event horizon in the larger image as it goes through 
the petal centers and tangency points. The petals are green and blue 
in the larger image, corresponding, as in \F{Shift2pi}, to which twin 
they are tangent to. The jump circles $j_1,j_2$ are also labeled.

\begin{figure}[h]
\begin{overpic}[width=.6\textwidth
%,grid,tics=10
]{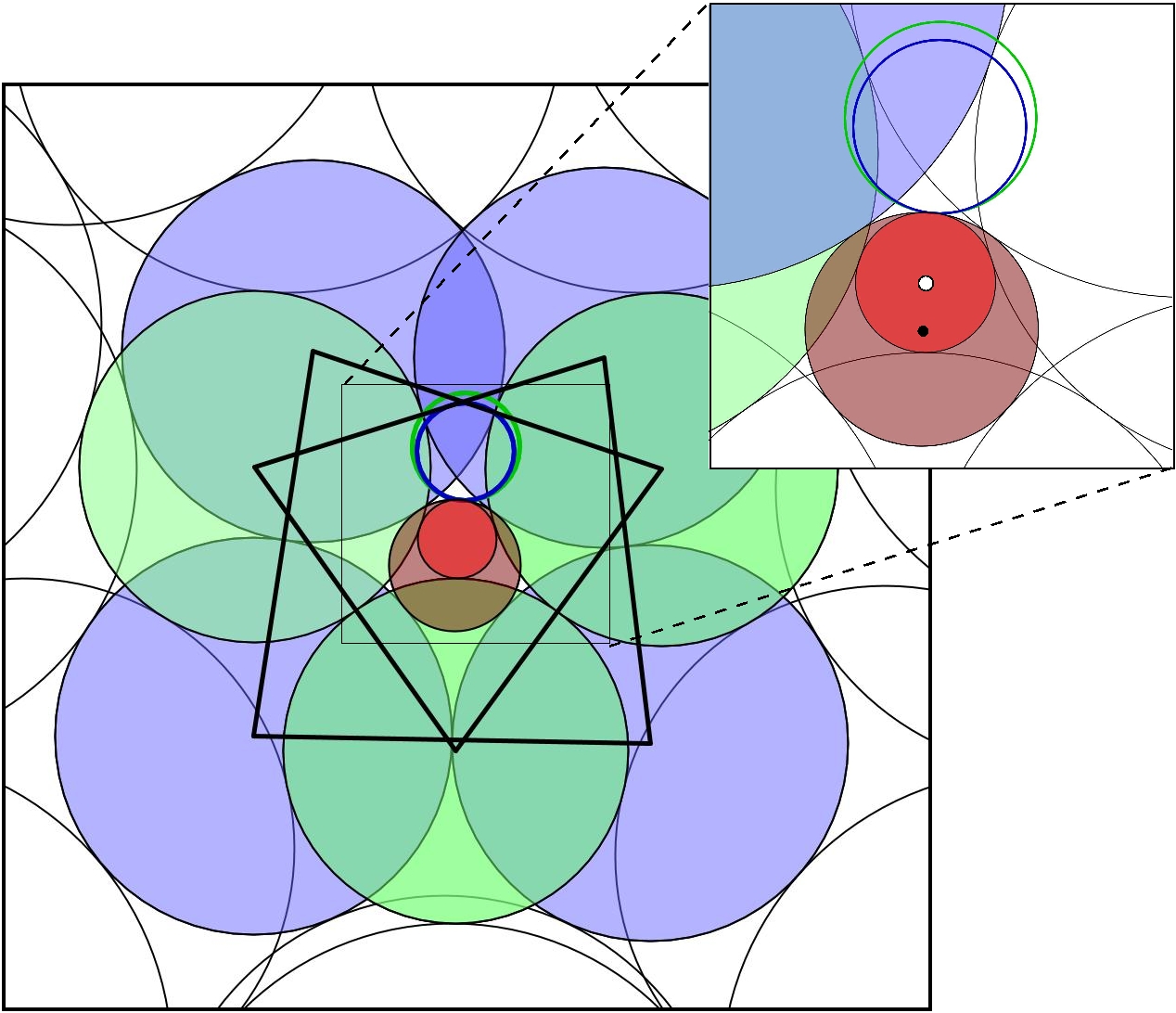}
\put (79,63) {$w$}
\put (9,43) {$j_1$}
\put (30,66) {$j_2$}
\end{overpic}
  \caption{Image circle packing in the neighborhood of 
shifted branching, with detail zoom.}
  \label{F:ShiftBranch}
\end{figure}

As with the singular branch image, the configuration of
\F{ShiftBranch} makes sense if one treats it as the projection of
circles lying on a two-sheeted surface $S$ branched over $w$. To see
this, consider the twins in the detail zoom: $t_1$ is the larger twin,
with center at the black dot and radius $r_1$. The smaller twin has
center at $w$ and radius $r_2<r_1$. Now imagine attaching a string of
length $r_1$ at the black dot and using it to draw the circle for
$t_1$ on $S$. As the string sweeps around, it will snag on the white
dot at $w$ and, like a yo-yo, trace out the smaller twin on $S$ before
finishing $t_1$. In other words, the union of the two twin circles
together is the projection of all points on $S$ which are distance
$r_1$ from the center of $t_1$ (that is, distance {\sl within}
$S$). Exactly this thought experiment was the genesis of shifted
branching.

If we disregard the chaperone circles and twins, the remaining circles
constitute a traditional tangency circle packing $P$ for $K$, with the
caveat that generically the circle for $v$ is ambiguous --- neither
the circle for $t_1$ nor for $t_2$ alone can serve as $c_v$. We need
to live with this ambiguity to achieve the branching behavior we want
outside the event horizon. (Having said this, there are (many)
settings which lead to identical twin circles, so $P$ then has this
common circle as $c_v$. All these configurations are identical and are
nothing but the circle packing we get when we choose traditional
branching at $v$.)

\vspace{10pt} This brings us to the matter of configuring black hole
combinatorics and parameters for this shifted branching; that is,
choosing the jump circles $j_1,j_2$ and their associated overlap
parameters $\gamma_1,\gamma_2$. We describe our scheme by referring to
\F{ShiftScheme}, which is the flower for $C_v$ in $P_K$. 

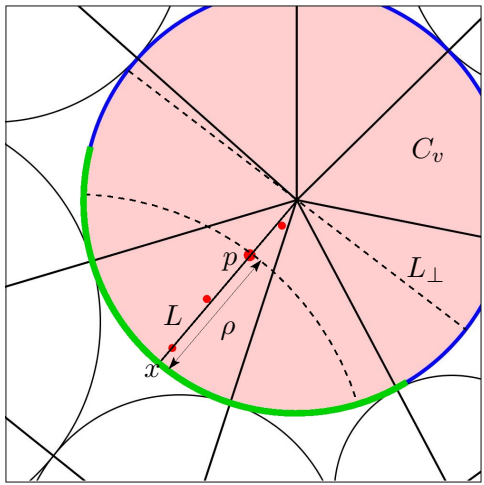
\begin{figure}[h]
\begin{overpic}[width=.5\textwidth
%,grid,tics=10
]{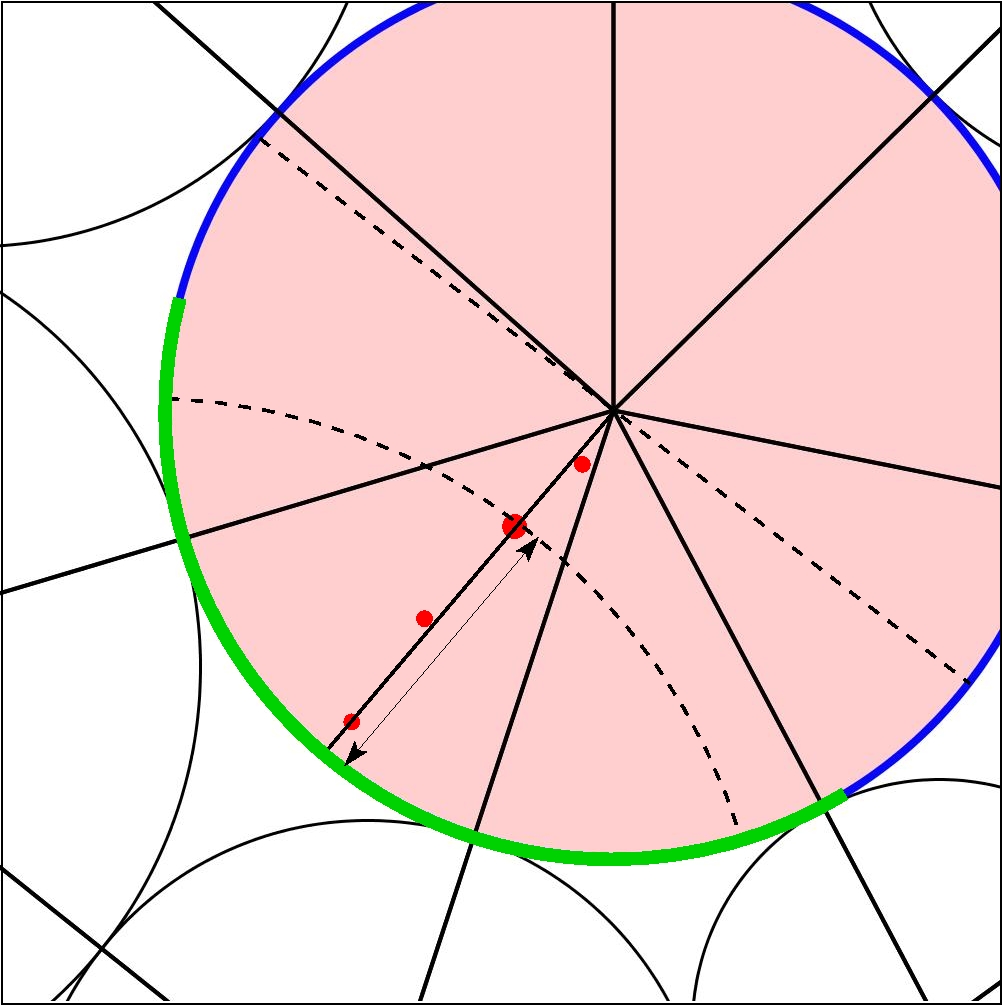}
\put (85,68) {$C_v$}
\put (29,22) {$x$}
\put (33,33) {$L$}
\put (84,43) {$L_{\perp}$}
\put (45.5,45.5) {$p$}
\put (45,31) {$\rho$}
\end{overpic}
  \caption{Choosing jumps and overlap parameters for shifted
    branching.}
  \label{F:ShiftScheme}
\end{figure}

The ultimate goal is to simulate branching at some point within $C_v$,
such as the indicated point $p$. In mapping to the branched image
packing $P$, the image of the boundary of $C_v$ wraps continuously
around the boundaries of both twin circles (as we described earlier in
referring to the branched surface $S$).  The jump circles and
parameters serve to split the boundary of $C_v$ into two arcs, the
blue one will be carried to $t_1$, the green, to $t_2$.

Here we need to observe how the jump and its parameter work together.
Recall that in the image packing, $\gamma_1\in[0,\pi]$ acts like a
dial: The value $\gamma_1=0$ forces $j_1$ to be tangent to both twins.
As $\gamma_1$ increases, it pushes $j_1$ away from $t_1$ and further
onto $t_2$.  When $\gamma_1$ reaches $\pi$, it forces the
counterclockwise petal $w_1$ to become tangent to $t_2$. This is a
transition point --- at this juncture, we could designate $w_1$ as the
jump circle and reset $\gamma_1$ to $0$ without altering anything in
the image packing. By then increasing the new $\gamma_1$ with the new
jump circle, we could push yet more boundary onto $t_2$.  In summary,
then, our circle packing map pushes more of $C_v$ onto $t_2$ by
increasing $\gamma_1$ and/or moving the designated jump $j_1$
clockwise.  Likewise, on the other side it pushes more of $C_v$ onto
$t_2$ by decreasing $\gamma_2$ and/or moving the designated jump $j_2$
counterclockwise.

To illustrate with the point $p$ of \F{ShiftScheme}, the scheme uses
the various labeled quantities: The point $x$ where the radial line
$L$ from the center of $C_v$ through $p$ hits $C_v$; the distance
$\rho$ from $p$ to $x$; the circular arc (dashed) through $p$ and
orthogonal to $C_v$; and the diameter $L_{\perp}$ perpendicular to
$L$.

To inform our choice of jumps and parameters, we take inspiration from
the properties of the branch value $w$ in the eventual image packing
--- that is, the center of the smaller twin, $t_2$. In qualitative
terms, the blue arc of $C_v$ should map to $t_2$, the rest of $C_v$ to
$t_1$. The point $x$ should map to the point of $t_2$ antipodal to the
tangency point of $t_1$ and $t_2$. The ratio of $\rho$ to the radius
of $C_v$ should reflect the ratio of the radii of the two twins. Thus,
when $p$ moves close to $C_v$, twin $t_2$ gets smaller, while as $p$
approaches the center of $C_v$, the radius of $t_1$ approaches that of
$t_2$. There is no way to ensure these outcomes precisely --- one
cannot know, {\sl a priori}, the outcomes in the image packing, as all
the circles get new sizes during computation. We will not burden the
reader with the messy details, but we have implemented methods which
realize these qualitative behaviors. We illustrate for $p$ and the
other three red branch points later, \F{ShiftMovie}.

\section{Fixing an Ahlfors Function}\label{S:FixAhlfors}
After successfully constructing a discrete Ahlfors function $\fw$ for
a combinatorial annulus $K$ in \S\ref{SS:DiscAhlfors}, we showed in
\S\ref{S:DiscIssue} how easily that construction can fail. Making
small modifications to $K$ that broke its translational symmetry, we
obtained a new combinatorial annulus $K'$ which does not support a
discrete Ahlfors function. The problem is non-trivial holonomy, and we
illustrated in \F{FailedAhlfors} with an attempt at traditional
branching using the same midline vertices $v_1,v_2$ we had used for
$\fw$.

It seems clear that for $K'$ the missing translational symmetry can be
blamed for the failure. We now apply the flexibility of generalized
branching to repair the damage.  Since $K'$ still has a midline and
reflective symmetry across it, we adopt the following strategy:
proceed with traditional branching at vertex $v_1$, but use shifted
branching near $v_2$. Symmetry simplifies our search for the correct
branching parameters in the black hole for $v_2$: namely, if we choose
vertices $j_1$ to be symmetric with $w_2$ across the midline, and
likewise, $j_2$ symmetric with $w_1$, and if we specify
$\gamma_2$=$\pi-\gamma_1$, then the shifted branch value must remain
on the midline. After some experimental tinkering, one can in fact
annihilate the holonomy and replicate the success we saw for the
original complex $K$ --- the process works. We do not show the image
packing $P$ because it is essentially indistinguishable from \F{AnnulusP}.
The point is that we are able to make the red cross-cuts coincident.

Admittedly, the fix was (almost) in for this example: we depended on
reflective symmetry to reduce the parameter search from a two- to a
one-dimensional problem. Nonetheless, it demonstrates well the need
and potential for generalized branching. We close by discussing the
broader issues.

\section{Parameter Space}\label{S:Parameter}
This paper is a preliminary report on work in progress. We have
focused on generalized branching at a single point $p$ in the interior
of $P_K$. The location of $p$ is continuously parameterized --- e.g, by
its $x$ and $y$ coordinates.  We have defined discrete generalized
branching which seems to handle patches of this parameter space.
Thus, when $p$ lies in an interior interstice, singular branching
involves two real parameters, $\gamma_1,\gamma_2$. When $p$ lies in an
interior circle, shifted branching involves jump circles and
parameters, but in our description of the mechanics it is clear that
this, too, is just two real parameters. The continuity of these
parameterizations may be phrased in terms of the branched packing
labels $R$ restricted to vertices on and outside of the event
horizon. 

While a proof remains elusive, experiments strongly suggest that this
continuity does hold. For example, \F{SingMovie} displays the branched
circle packings associated with branching at the four red dots in
\F{SingScheme}, progressing from lower left to upper right (the
third of these is the packing for the distinguished point $p$ from
\F{SingScheme}). The branch value is roughly at the center in each
image. Subject to this and related normalizations, the radii and
centers of $P$ appear to be continuous in $\gamma_1,\gamma_2$.

\begin{figure}[h]
\begin{overpic}[width=1.0\textwidth
%,grid,tics=10
]{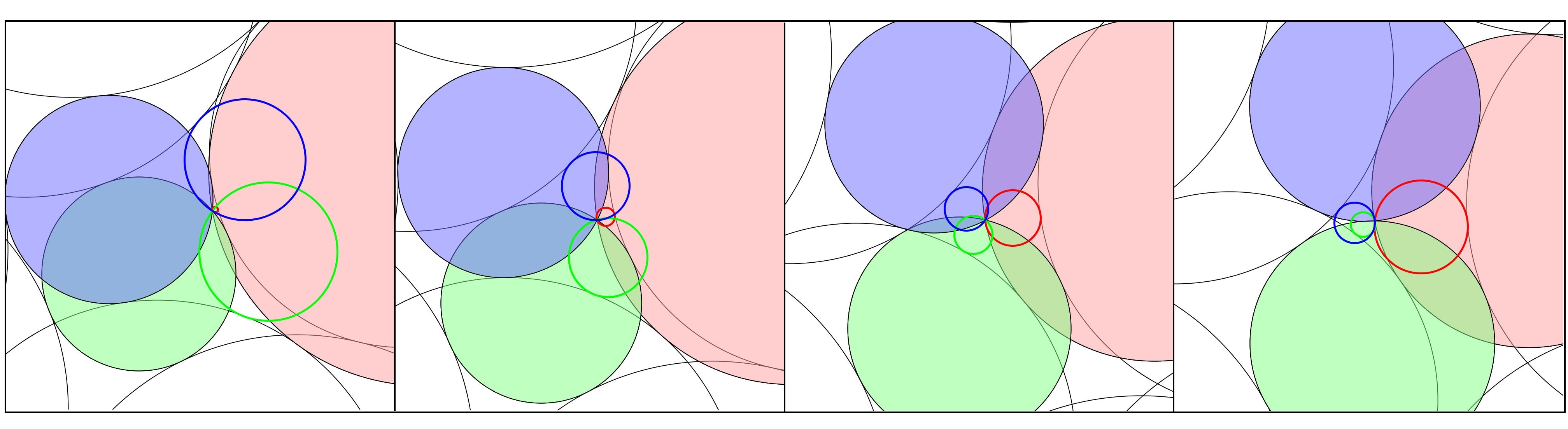}
%\put (2,2) {(a)}
\end{overpic}
  \caption{Singular branching for the four
red branch points of \F{SingScheme}.}
  \label{F:SingMovie}
\end{figure}

\F{ShiftMovie} provides a similar sequence of shifted branched
packings for the four red dots of \F{SingScheme} (caution: the 
chaperones play different roles now). Again we have positioned
the branch values roughly at the center in each image; the third
one corresponds to \F{ShiftBranch}. Here, too,
experiments suggest continuity in radii and centers as we manipulate
the two shifted branching parameters.

\begin{figure}[h]
\begin{overpic}[width=1.0\textwidth
%,grid,tics=10
]{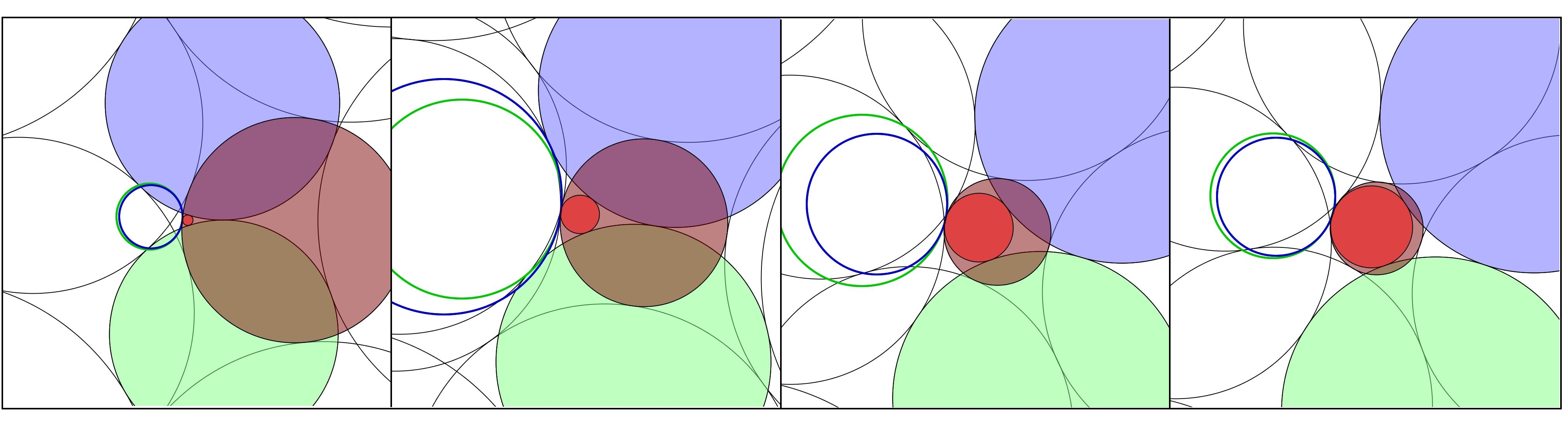}
%\put (2,2) {(a)}
\end{overpic}
  \caption{Shifted branching for the four
    red branch points of \F{ShiftScheme}.}
  \label{F:ShiftMovie}
\end{figure}

Concatenating the 8 frames in these last two figures highlights
another parameterization issue: How are our various patches of
parameter space sewn together? If $p$ lies on the mutual boundary of a
circle and an interstice, for instance, its generalized branched
packing may be treated as a limit of either singular branching from
the interstice side or shifted branching from the circle side.  We
have {\sl ad hoc} methods for such transitions, though we have yet to
formalize the details of parameter alignment. Nevertheless, our images
may give a feel for the transition: The interstice formed by
$\{C_{v_1},C_{v_2},C_{v_3}\}$ in \F{SingScheme} is contiguous to the
circle $C_v$ of \F{ShiftScheme}; that is, $v=v_1$. So the 8 frames
from \F{SingMovie} and \F{ShiftMovie} together are part of a movie as
the branch point transitions from singular to shifted.  Image circles
$\{c_{v_1},c_{v_2},c_{v_3}\}$ remain red, green, and blue,
respectively, throughout these 8 frames. In the last frame from
\F{SingMovie} note that these three appear to be in clockwise order
(as we discussed earlier).  Compare this to the first frame of
\F{ShiftMovie}: the red circle has now split into twins, with the
branch value in the smaller twin, so the (small) red, green, and blue
are again correctly oriented --- the branch point has successfully
punched though from the interstice to the circle, and in the last
frame of \F{ShiftMovie}, it is nearing traditional branching at
$v$. This is the type of experimental evidence supporting our
contention that the two parameter patches can be aligned to maintain
continuity.

\section{Global Considerations}\label{S:Global}
We stated in the introduction that our aim is to bridge the principal
gap remaining in discrete function theory, namely the existence and
uniqueness of discrete meromorphic functions. Although we have local
machinery, we have not confronted the global problem head-on.  A few
words are in order.

Naturally, one of the first goals would be a more complete theory for
discrete rational functions, branched mappings from $\bP$ to
itself. Here $K$ would be a combinatorial sphere and one would need
$2n$ branch points for a mapping of valence $n+1$. There is a
tantalizing approach based on Oded Schramm's metric packing theorem
which has motivated some of our work.  In \cite{oS91a}, Schramm proves
remarkable existence and uniqueness results for packings of ``blunt
disklike'' sets --- for instance, circles defined using a Riemannian
metric on $\bP$.  Suppose, then, that we are given a classical
rational function $\fF:\bP\rightarrow \bP$. We can define a metric $d$
on $\bP$ as the pullback under $\fF$ of the spherical metric on
$\bP$. Finding a packing for $K$ by circles in this metric $d$ is
tantamount to finding a normal circle packing $P$ on $\bP$, and the
map from $P_K$ to this $P$ would be our discrete rational
function. Unfortunately, at the critical values of $\fF$ the pullback
metric $d$ is not Riemannian; the direct analogue of Schramm's result
does not hold, as can be seen, for example, with circles that
degenerate. There is still some hope, however, as our constructions
demonstrate --- the twin circles of a shifted branch point are, after
all, the image of a single circle in a pullback metric $d$.

Our hands-on approach still faces many hurdles in practice. On the
sphere, for instance, there is no packing algorithm --- Perron methods
rely on the monotonicity of Lemma~\ref{L:Ashe}, which fails in the
positive curvature setting.  And in other settings, such as Ahlfors
and Weierstrass, we have had to depend on symmetry. A generic
combinatorial torus $K$ is likely to have no Weierstrass function
using traditional branching. Though we believe generalized branching
provides the flexibility to overcome the holonomy obstructions, early
attempts have faltered due to the curse of (even small) dimension: we
don't yet know how to search a two-dimensional space for parameters that
will annihilate non-trivial holonomies. We succeeded in the Ahlfors
case because partial symmetry reduced us to a one-dimensional search.

We face other global difficulties as well. We list a few. We have
restricted attention to {\sl simple} branching; at least in the case
of shifted branching, one can see a chance to allow higher order
branching --- replacing twins with triplets, etc. In general, one also
needs to allow branching at more than one point, but the existence of
branched packings then encounters global combinatorial issues. The
notion of black holes will also need to be extended, since
combinatorics may lead to patches of degenerate radii (versus isolated
degenerate radii) for branch points in certain combinatorial
environments.

In other words, there is considerable work to be done. Nevertheless,
we contend that discrete generalized branching addresses --- in theory
if not in practice --- the key obstruction remaining in discrete
analytic function theory. This obstruction, of course, is not the only
one --- so get to work, David.

\bibliographystyle{amsplain}
\bibliography{Gen-bib}

\end{document}